\documentclass[11pt]{amsart}
\usepackage{amsmath}
\usepackage{amssymb}
\usepackage{euscript}
\usepackage{pifont}
\usepackage{tikz}
\usepackage[all]{xy}
\numberwithin{equation}{section}
\numberwithin{figure}{section}

\input cyracc.def

\usepackage{amssymb,amsmath}
\usepackage{pstricks}
\usepackage{eepic}

\newrgbcolor{blue}{0 0 .7}

\newtheorem{definition}[equation]{Definition}
\newtheorem{theorem}[equation]{Theorem}
\newtheorem{proposition}[equation]{Proposition}

\newtheorem{lemma}[equation]{Lemma}
\newtheorem{example}[equation]{Example}
\newtheorem{remark}[equation]{Remark}

\newcommand{\w}{\wedge}
\newcommand{\n}{\notag}
\newcommand{\noi}{\noindent}

\newcommand{\hook}{\hookrightarrow}
\newcommand{\tb}{\textbf}
\newcommand{\vsp}{\vspace}
\newcommand{\sq}{ \;  \square}

\newcommand{\C}{ \mathbb{C} }

\newcommand{\W}{ \mathcal{W} }

\newcommand{\PP}{ \mathbb{P} }

\newcommand{\K} { \mathcal{K} }

\newcommand{\del}{\delta \phi}

\newcommand{\Pf} {\noi \emph{Proof.}$\; \;$}

\newcommand{\beq}{ \begin{equation}}
\newcommand{\eeq}{ \end{equation} }

\begin{document}

\sloppy
\thispagestyle{empty}

\title{On Gronwall conjecture}
\author{Joe S. Wang}
\address{St. Louis, MO 63124  USA}
\email{jswang12@gmail.com}
\begin{abstract}
Gronwall conjecture states that
a  planar 3-web
which admits more than one distinct linearization
is locally equivalent to an algebraic web.
We give a partial answer to the conjecture in the affirmative
for the class of  planar 3-webs
with the web curvature that vanishes to order three at a point.
The differential relation on the third order jet of web curvature
provides  an explicit criterion  for unique linearization.
\end{abstract}
\maketitle

\tableofcontents

%-------------------------------------------------
%-------------------------------------------------
\section{Introduction}\label{sec1}
A planar $\, d$-web is by definition
a set of $\, d$ transversal foliations by curves
on a two dimensional surface.
Let $\, \W$ be a $\, d$-web on a surface $\, M$.
Let $\, \PP^2$ be the projective plane.
A \emph{linearization} of $\, \W$
is an immersion $\, M \hook \PP^2$ such that
each leaf of the foliations is mapped to a line.
Two linearizations are equivalent if they are isomorphic up to
projective transformation of $\, \PP^2$, and otherwise distinct.
A web on an open subset of $\, \PP^2$ is linear
when each leaf of the foliations is a part of a line.

There exist  a distinguished class of  linear  webs.
Let $\, \Gamma \subset  (\PP^2)^*$ be a reduced, degree $\, d$
algebraic curve in the dual projective plane.
On a neighborhood of a generic point of $\, \PP^2$,
$\, \Gamma$ induces a $\, d$-web $\, \W_{\Gamma}$ by the standard dual construction.
Such a linear web on an open subset of $\, \PP^2$ is called \emph{algebraic}.
It is well known that a planar 3-web is locally equivalent to
an algebraic web when its web curvature vanishes, \cite{Ch2}\cite{He}.

\vsp{.5pc}\noi
\tb{Gronwall conjecture.}\;\;
Let $\, \W$ be a planar 3-web
which admits more than one
distinct  linearization.
Then $\, \W$ is locally equivalent to an algebraic web.

\vsp{.5pc}
\noi
Since a plane cubic curve   has local invariants,
an algebraic 3-web does admit more than one distinct linearization.

The purpose of  the present paper
is to give a partial answer to the conjecture in the affirmative.
We  provide  an  explicit criterion for unique linearization
in terms of a differential relation
on the third order jet of web curvature.

\vsp{.5pc}\noi
\tb{Main theorem.}\;\;
Let $\, \W$ be a planar 3-web on a connected surface $\, M$
which admits two distinct linearizations.
Let   x$_0\in M$ be a reference point.
Suppose the web curvature of $\, \W$ vanishes to order three at x$_0$.
Then the web curvature vanishes identically,
and
$\, \W$ is locally equivalent to an algebraic web.

\vsp{.5pc}
\noi
If a planar 3-web $\, \W$  has the web curvature
which vanishes at least  to order three at a point
but  which does not vanish identically,
then $\, \W$ admits at most one distinct local  linearization.
We give a construction of
examples of such linear 3-webs
with the web curvature that vanishes to arbitrary order at a point.

Let x$_0 \in \PP^2$ be a reference point.
Let $\, \mathfrak{M}$ be the moduli space of germs of linear 3-webs at x$_0$.
Since a germ of linear foliation at x$_0$ is determined by a germ of curve in $\, (\PP^2)^*$,
there are roughly three arbitrary functions of one variable worth linear 3-webs in $\, \mathfrak{M}$.
Let $\, K$ denote the web curvature function on  $\, \mathfrak{M}$
(web curvature is a relative invariant of a 3-web, \eqref{curvature}.
One may determine a section
by a choice of germ of nonzero 2-form at x$_0$).
Consider the valuation map
\beq
j^3 K:   \mathfrak{M} \to \C^{10}, \n
\eeq
that records the third order jet of the web curvature of the linear 3-web at x$_0$.
%It is easily checked that this map is surjective.
Main theorem  answers Gronwall conjecture in the affirmative
for the  subset of linear 3-webs $\, (j^3 K)^{-1}(0)$
of codimension at most 10  in  $\, \mathfrak{M}$.
The  differential relation satisfied by the web curvature, \eqref{K3functional},
may provide a basis for further in depth
analysis toward  Gronwall conjecture.

In the article  \emph{'Sur les equations entre trois variables representables
par les nomogrammes  a points aligne'} published  in 1912,
Gronwall considered the following problem, \cite{Gro}.
Let $\,\{ x_1, \, x_2, \, x_3 \}$ be a set of three functions on a surface $\, M$
that satisfy a relation
\beq\label{Frelation}
R(x_1, \, x_2, \, x_3)=0.
\eeq
When does there exist
 three colinear curves
$\, \Gamma^i(x_i) : M \to (\PP^2)^*, \, i=1, \, 2, \, 3 \,$?
Such a set of three curves is called \emph{a nomographic representation} of  \eqref{Frelation}.
Note that the associated map  $\, M \to \PP^2 = ((\PP^2)^*)^*$
gives  a linearization of the 3-web
defined by   $\,\{ x_1, \, x_2, \, x_3 \}$.

Let $\,  \Gamma^i(x_i) = ( f_i(x_i), \, g_i(x_i), \, 1 )$ be a parametrization
in an affine chart.
$\, \Gamma^i$'s are colinear when there exist  two functions $\{ u, \, v \}$ on $\, M$
such that
\beq\label{nomograph}
g_i(x_i) = u \, f_i(x_i) + v, \; \;  i = 1, \, 2, \, 3.
\eeq
By successive differentiation of  \eqref{nomograph},
one can eliminate $\, f_i(x_i), \, g_i(x_i)$
and derive   two fourth order PDE's for  $\{ u, \, v \}$.
Gronwall showed, among other things, that solvability of this pair of   PDE's
is also a sufficient condition to admit a nomographic representation.

On page 61 of  \cite{Gro}, he wrote;
\emph{
Dans un travail ulterieur,
je formerai explicitement
l'integrale commune des equations aux derivees partielles
du paragraphe 1,
et je ferai voir que le cas du
paragraphe 4 est le seul ou
l'equation donnee admette les
representations nomographiques
essentiellement distinctes.
}

Here
\emph{'des equations aux derivees partielles du paragraphe 1' }
means the aforementioned pair of fourth order PDE's,
and  \emph{ 'le cas du paragraphe 4'}
means the case when each $\, \Gamma^i(x_i), \, i = 1, \, 2, \, 3$, is linear.
The alluded subsequent work does not appear to be published.

The compatibility equations to admit a linearization
impose  a stringent set of conditions on a 3-web.
Bol, and Boruvka showed that
a  planar 3-web with nonzero web curvature
admits at most 16 distinct linearizations, \cite{Bol0}\cite{Bor}.
Vaona improved this bound to 11, \cite{Vao}.
Smirnov gave a proof of Gronwall conjecture,
the content of which is not available to us, \cite{Smi}\footnotemark.
\footnotetext{
We  were unable to  locate the original paper
\cite{Smi}  in any format.
A review  is available at   Zbl 0261.53007.}

More recently,
Grifone,  Muzsnay, and Saab   proved the bound of 15 distinct linearizations, \cite{GMS}.
Goldberg and Lychagin also  proved the bound of 15
in relation to their work on Blaschke conjecture for 3-webs, \cite{GL}.
These results were obtained essentially by
determining a bound on the number of common roots
of  a  set of polynomial compatibility equations
through evolved, intricate differential algebraic analysis.
They both assumed  that the planar 3-web
has nonzero web curvature.

Let us give the outline of  proof  of   Main theorem.
Let   $\, \PP(TM) \to M$ be the projective tangent bundle of a surface $\, M$.
Let $\, \mathcal{D}$ be the canonical contact 2-plane field on $\, \PP(TM)$.
A \emph{path geometry}  on   $\, M$
is by definition
a $\, \mathcal{D}$-horizontal foliation on $\, \PP(TM)$
transversal to the fibers of  the projection $\, \PP(TM) \to M$.
Under the projection  to $\, M$,
it determines a unique path, or a curve,
tangent to each direction $\, [\mbox{v}] \in \PP(TM)$.
A path geometry is flat when it is locally equivalent to
the standard path geometry of  the projective plane.

Let $\, \W$ be a $\, d$-web on $\, M$.
A linearization of $\, \W$ up to projective transformation
is equivalent to
a flat path geometry on $\, M$ up to isomorphism
such that
each leaf of the foliations of $\, \W$ is a part of a path.

The fundamental observation for our investigation is that
two distinct flat  path geometries intersect
along a generalized 3-web
defined by the base locus of a conformal class of symmetric cubic differential.
This means that if a 3-web $\, \W$ admits two distinct linearizations,
the two linearizations in turn uniquely determine $\, \W$.
The condition that $\, \W$ is linear imposes additional set of compatibility equations,
which allow  one to close up the structure equation  for $\, \W$
by the over-determined PDE machinery.

An examination of the resulting structure equation
reveals the following rigid property of the web curvature.
Let $\, K$ denote the web curvature of $\, \W$.
Let $\, \mathcal{K}^s = \langle \, K, \, \nabla K, \, \nabla^2 K, \, ... \,  \nabla^s K \, \rangle$,
$\, s = 0, \, 1, \, 2, \, ... \, $,
be the ideal of functions generated by
the $\, s\,$th-jet of $\, K$
(by definition of the relative invariant $\, K$,
the ideal $\, \mathcal{K}^s$ is well defined).
The structure equation implies the differential relation
\beq\label{K3functional}
d \mathcal{K}^3 \equiv 0, \mod  \; \mathcal{K}^3.
\eeq
Main theorem  follows  by the uniqueness theorem of ODE.

In Section \ref{sec11},
we give a definition, and list the basic properties of   path geometry structure on a surface.
In Section \ref{sec2},
an analysis of Maurer-Cartan equation for the deformation of  flat path geometry
on an open subset of $\, \PP^2$
yields a conformal class of symmetric cubic differential $\, \sigma$, \eqref{sigma}.
A global consideration using $\, \sigma$ shows that
there exists a unique flat path geometry on $\, \PP^2$, Theorem \ref{P2}.
In Section \ref{sec3},
we impose the condition that the 3-web $\, \W_{\sigma} = \sigma^{-1}(0)$ is linear.
An explicit formula for the web curvature of  $\, \W_{\sigma}$ is obtained
as a fifth order invariant of the deformation, \eqref{webcurvature}.
Moreover,
the structure equation for  $\, \W_{\sigma}$ closes up
at order eight.
In Section \ref{sec4},
a direct computation using
the closed structure equation for $\, \W_{\sigma}$ implies that
the curvature ideal $\, \mathcal{K}^3$ is differentially closed.
When $\, \W_{\sigma}$ contains one, or two pencils,
the curvature ideal $\, \mathcal{K}^2$,  or $\, \mathcal{K}^1$
is differentially closed respectively, Theorem \ref{main}.
In Section \ref{sec5},
we give a remark toward the full proof of Gronwall conjecture.

For  reference on path geometry on a surface,
we cite  \cite{Br}\cite{BGH}.
For general reference on web geometry,
we cite \cite{Ch2}\cite{GS}\cite{PP}.
On the linearization of planar webs,
we cite \cite{He}\cite{GMS}\cite{GL}, and the references therein.

The method of moving frames, and  exterior differential systems are used
throughout the paper without specific reference.
For the standard reference,
we cite \cite{Ga}\cite{BCG}\cite{IL}.

For a uniform treatment, we adopt the complex, holomorphic category.
All of the results are valid in the  real, smooth category with minor modifications.

The majority of computations were performed using
the computer algebra system \texttt{Maple} with \texttt{difforms} package.

\subsection{Path geometry}\label{sec11}
Let $\, M$ be a two dimensional manifold.
Let $\, \PP(TM) \to M$ be the projective tangent bundle
equipped with the canonical contact 2-plane field
$\, \mathcal{D} \subset T(\PP(TM))$, \cite{Br}.
\begin{definition}
A \emph{path geometry}  on a surface $\, M$ is
a $\, \mathcal{D}$-horizontal foliation
transversal to the fibers of  projection $\, \PP(TM) \to M$.
\end{definition}

Let $\, (x, \, y)$ be a generic local coordinate of $\, M$.
Introduce a variable $\, p\,$ so that
$\, (x, \, y, \, p)$ is a local coordinate of $\, \PP(TM)$,
and that
$\, \mathcal{D} = \langle \, dy - p\, dx \, \rangle^{\perp}$.
Let $\, \mathcal{F}$ be a $\, \mathcal{D}$-horizontal foliation
that defines a path geometry.
By transversality condition,
there exists a function $\, f(x, \, y, \, p)$ such that
$\, \mathcal{F}$ is locally defined by the corank one Pfaffian system
\beq
\langle \, dy - p\, dx, \, \, dp - f(x, \, y, \, p)\, d x \rangle. \n
\eeq
The paths of the path geometry $\, \mathcal{F}$
are locally the integral curves of the second order ODE
\beq\label{definingODE}
\frac{d^2 y}{d x^2} = f(x, \, y, \, \frac{dy}{dx}).
\eeq

The flat model of path geometry is
the standard  homogeneous  path geometry
of lines on the projective plane.
Let  $\, SL_3 \C$ be the group of 3-by-3 matrices of determinant one.
Let $\, P_1, \, P_2$, and $\, P$ be    the following   subgroups.
\begin{align}
P_1 &= \{ \,
\begin{pmatrix}
* & * & * \\
\cdot & * & * \\
\cdot & * & *
\end{pmatrix}
\, \}, \n \\
P_2 &= \{ \,
\begin{pmatrix}
* & * & * \\
* & * & * \\
\cdot & \cdot & *
\end{pmatrix}
\, \}, \n \\
P& = P_1 \cap P_2, \n
\end{align}
where '$\, \cdot$' denotes 0 and '$\, * $' is arbitrary.
Consider the double fibration.

\begin{picture}(300,87)(-67,-14)
\put(147,40){$\PP(T\PP^2) = SL_3 \C / P$}
\put(175,25){$\searrow$}
\put(145,25){$\swarrow$}
\put(184,10){$ (\PP^2)^* =SL_3 \C / P_2$ }
\put(71,10){ $SL_3 \C / P_1 = \PP^2 $ }
\put(187,30){$\pi_2$}
\put(133,30){$\pi_1$}
\end{picture}

\noi
The foliation by fibers of $\, \pi_2$ induces the standard     path geometry on $\, \PP^2$,
which is locally described by the second order ODE
\beq
\frac{d^2 y}{d x^2} =0. \n
\eeq
The foliation by fibers of $\, \pi_1$ induces the dual path geometry on $\, (\PP^2)^*$.

Cartan, led by his geometric study of differential equations,
solved the local equivalence problem for path geometry, \cite{Ca}\cite{BGH}.
%Let $\, \mathcal{F}$ be the foliation on $\, \PP(TM) \to M$ that defines a path geometry on $\, M$.
Let $\, \mathcal{F}^*$ be the foliation by the fibers of the projection $\, \PP(TM) \to M$.
Assume the moduli space of leaves $\, \PP(TM) / \mathcal{F}$ is a smooth manifold,
and
consider the following incidence double fibration.

\begin{picture}(300,87)(-67,-14)
\put(147,40){$\PP(TM)$}
\put(175,25){$\searrow$}
\put(145,25){$\swarrow$}
\put(184,10){$ M^* =\PP(TM) / \mathcal{F}$ }
\put(61,10){ $\PP(TM) / \mathcal{F}^* = M$ }
\put(187,30){$\pi_2$}
\put(133,30){$\pi_1$}
\end{picture}

The local equivalence problem   is solved on
a $\, P$-bundle $\, B \to \PP(TM)$
together with an $\, \mathfrak{sl}_3 \C$-valued Cartan connection form $\, \phi$.
The bundle $\, B \to \PP(TM)$ can be considered as the union of infinitesimal homogeneous spaces
$\,  SL_3 \C \to SL_3 \C / P$ connected by Cartan connection $\, \phi$.

The following theorem is drawn from \cite[p176]{Br}.
Let $\, (\phi_{i,j})_{i,j=0}^2 $ denote the components of $\, \phi$.
Let $\, \mathfrak{p} \subset  \mathfrak{sl}_3 \C$ be
 the Lie algebra of the  subgroup $\, P \subset SL_3 \C$.

\begin{theorem}[Cartan]\label{defining}
Let $\, \mathcal{F}$ be the foliation on $\, \PP(TM) \to M$
that defines a path geometry.
There exists a principal right $\,P$-bundle $\, \tau: B \to   \PP(TM) $
and an $\, \mathfrak{sl}_3 \C$-valued  1-form  $\, \phi$
on $\,B$ with the following  properties:

(1) For each $\, b \in B$, the map $\, \phi_b: T_b B \to  \mathfrak{sl}_3 \C$
is an isomorphism and
$\, \phi$ pulls back to each fiber of  $\, \tau$
to be the canonical $\, \mathfrak{p}$-valued left-invariant 1-form.

(2) $\, R_g^* \, \phi = g^{-1} \, \phi \, g$ for each $\, g \in P$.
Here $\, R_g$ is the right action by $\, g$.

(3) For some (and hence any) section $\,\textnormal{u}: \PP(TM) \to B$,
the pullback 1-form $\, \varphi =  \textnormal{u}^* \, \phi$
has the properties that
the leaves of the foliation $\, \mathcal{F}$
are the integral curves of  $\, \varphi_{2,0} = \varphi_{2,1} = 0$
while the 1-form $\, \varphi_{1,0}$ is nonzero on each leaf,
and that
the leaves of the foliation $\, \mathcal{F}^*$
are the integral curves of  $\, \varphi_{1,0} = \varphi_{2,0} = 0$
while the 1-form $\, \varphi_{2,1}$ is nonzero on each leaf.

(4) The curvature 2-form   $\, \Phi = d \phi + \phi \w \phi$ satisfies
\beq
\Phi =
\begin{pmatrix}
\cdot &  L_1  \, \phi_{1,0} \w \phi_{2,0}  & \Phi_{0,2}  \\
\cdot  & \cdot                                           &  L_2 \, \phi_{2,1}  \w \phi_{2,0}  \\
\cdot  & \cdot  & \cdot
\end{pmatrix}\n
\eeq
for some functions $\, L_1$ and $\, L_2$  on $\, B$.
$\, \Phi$ vanishes if and only if the path geometry is locally equivalent
to the flat model on the projective plane.

The pair
$\, (B, \, \phi)$ is uniquely characterized by these four properties:
If    $\, (B', \, \phi')$
also satisfies them
then there exists a unique bundle isomorphism
$ E : B \to B'$ covering the identity on $\, \PP(TM)$
so that  $ E^*\phi' = \phi$.
\end{theorem}

In terms of the second order ODE \eqref{definingODE},
the invariant $\, L_2$ vanishes when $\, f ( x, \, y, \, p)$ is at most cubic in $\, p$.
In this case, \eqref{definingODE} is the equation of geodesics of a projective connection on $\, M$.
Dually,
the invariant $\, L_1$ vanishes when $\, f ( x, \, y, \, p)$ satisfies
\beq\label{Mvanish}
\frac{d^2}{dx^2}f_{pp} -4 \, \frac{d}{dx} f_{py}
+f_p( 4\, f_{py} -  \frac{d}{dx}f_{pp})
-3\, f_y \, f_{pp} + 6\, f_{yy} =0.
\eeq
In this case, the dual equation of \eqref{definingODE}
 is the equation of geodesics of a projective connection on $\, M^*$.
By Bianchi identity
\beq
d\Phi + \phi \w \Phi - \Phi \w \phi = 0, \n
\eeq
$\, \Phi$ vanishes when $\, L_1$ and $\, L_2$ vanish.

Let $\, \W$ be a $\, d$-web on a surface $\, M$.
Let $\, (x, y)$ be a generic local coordinate of $\, M$ such that
$\, \W$ is  defined by $\, d$ first order ODE's
\beq
\frac{dy}{dx} = p^a(x, y), \; \, a = 1, \, 2, \, ... \, d, \n
\eeq
where $\, p^a \ne p^b$ for $\, a \ne b$.
The analysis above implies that
a local linearization of $\, \W$ is equivalent to finding
a second order ODE of the form
\begin{align}\label{fiscubic}
\frac{d^2y}{dx^2} &= f ( x, \, y, \,\frac{dy}{dx})   \\
&= h_3(x,y) \, (\frac{dy}{dx})^3 + h_2(x,y) \, (\frac{dy}{dx})^2
 +h_1(x,y) \, \frac{dy}{dx}  + h_0(x,y), \n
\end{align}
which satisfies  \eqref{Mvanish}, and
\beq\label{vanderm}
\frac{d}{dx} p^a(x, y) =h_3(x,y) \, p^a(x, y) ^3 +h_2(x,y) \, p^a(x, y) ^2
 +h_1(x,y) \, p^a(x, y)  + h_0(x,y), \; \; a=1, \, 2, \, ... \, d.
\eeq

When $\, d \geq 4$, \eqref{vanderm} and Vandermonde identity imply  that
$\, (h_3,   h_2,   h_1,   h_0)$ is uniquely determined by $\, p^a$'s.
Gronwall conjecture claims that when $\, d=3$ and
the web curvature of $\, \W$ is not identically zero,
\eqref{Mvanish}, \eqref{fiscubic}, and \eqref{vanderm}
 uniquely determine  $\, (h_3,   h_2,   h_1,   h_0)$.

%-------------------------------------------------
%-------------------------------------------------
\section{Deformation of flat  path  geometry}\label{sec2}
In this section,
we establish the fundamental  structure equation
for the deformation of flat path geometry
on an open subset of  the projective plane.
A  cubic differential   arises,
which encodes the essential local  information of deformation.
The  cubic differential is the unifying theme
throughout the paper.
It also provides the practical  computational perspective
for the  analysis in the later sections.

Let $\, V = \C^3$ be the three dimensional complex vector space.
Let $\, \PP^2 = \PP(V)$ be the projective plane.
Let $\, \PP(T\PP^2) \to \PP^2$ be the projective tangent bundle
equipped with the standard $\, SL_3 \C$-invariant  flat path geometry.
There exists a $\, P$-bundle   $\,  F= SL_3 \C \to  \PP(T\PP^2)$
with the $\, \mathfrak{sl}_3 \C$-valued Cartan connection form $\, \phi$,
which in this case is the left invariant Maurer-Cartan form of   $\, SL_3 \C$.
The pair $\, (F, \, \phi)$ satisfies
the defining properties  described in Theorem \ref{defining}.

Suppose $\,  \PP^2$ is given another flat path geometry structure.
Let $\, F' =  SL_3 \C  \to  \PP(T\PP^2)$ be the associated
$\, P'( \simeq P)$-bundle with Cartan connection form $\, \pi$.
From the definition of the bundle $\, F'$ and $\, F$,
one may regard  $\, F'$ as a graph over  $\, F$,
and assume $\, \pi$ is
another $\, \mathfrak{sl}_3 \C$-valued Cartan connection form on $\, F$.
In effect,
one may adopt the  following analysis
as the constructive definition of  $\, (F', \, \pi)$.

Set
\beq
\pi = \phi + \del. \n
\eeq
Maurer-Cartan equation  for $\, \pi$ and $\, \phi$  implies
the fundamental structure equation for the deformation $\, \del$;
\beq\label{fundast}
d(\del) + \del \w \phi + \phi \w \del + \del \w \del = 0.
\eeq
Differentiating the components of $\, \del$ from now on
would mean
applying this structure equation.

Set
\begin{align}
\phi_{2,0}&=\theta, \n \\
\phi_{1,0}&=\omega, \n \\
\phi_{2,1}&=\eta.\n
\end{align}
From the defining properties of Cartan connection form,
and using the group action by $\, P'$ on $\, \pi$,
one may assume
\begin{align}\label{No1}
\del_{2,0} &=0, \\
\del_{1,0} &=0,  \n \\
\del_{2,1} &\equiv 0, \mod \; \eta, \,  \omega. \n
\end{align}
Since $\,  d\theta  + \eta \w \omega \equiv 0,  \mod \; \theta$,
this forces
\beq\label{No2}
\del_{2,1}  \equiv 0,  \mod \;   \omega.
\eeq

At  this stage, the lower left corner of $\, \del$ is normalized such that
\beq\label{initial}
\del =
\begin{pmatrix}
* &    * & * \\
\cdot & *& * \\
\cdot  &  t \, \omega & *
\end{pmatrix},
\eeq
for a coefficient $\, t$.
The analysis in this section  will show that
the remaining coefficients of $\, \del $ are in the linear span of
the deformation function $\, t$ and its successive derivatives.
Two flat path geometries are  isomorphic, or $\, \pi = \phi$,
when $\, t$ vanishes identically.

The subgroup of elements of $\, P'$ preserving  \eqref{No1}, \eqref{No2}
are of the form (modulo  the finite subgroup of   center $\, Z(SL_3 \C) \simeq \mathbb{Z}_3$)
\beq
\begin{pmatrix}
1 &    \cdot & * \\
\cdot & 1& \cdot  \\
\cdot  & \cdot & 1
\end{pmatrix}. \n
\eeq
Using this group action, one may further normalize $\, \del\,$ so that
\beq\label{del00}
\del_{0,0} - \del_{2,2}  \;\;\; \mbox{has no $\theta$ component}.
\eeq
Equations \eqref{No1}, \eqref{No2}, and  \eqref{del00} in turn
uniquely determine $\, \del$.

We wish to determine the rest of the components of $\, \del$
by successive application of   the fundamental structure equation  \eqref{fundast}
to  the initial  state \eqref{initial}, \eqref{del00}.

\texttt{Step 1.}
Differentiating $\, \del_{2,0} = 0$, one gets
\beq
(\del_{0,0} - \del_{2,2} ) \w \theta=0. \n
\eeq
By \eqref{del00}, this implies
\beq\label{step1}
\del_{0,0} = \del_{2,2}. \n
\eeq
Since Trace($\del)=0$,  one has $\, \del_{1,1}=-2 \,  \del_{0,0}$.

\texttt{Step 2.} Differentiating $\, \del_{1,0} = 0,  \;   \del_{2,1} = t \omega$,
and after applying Cartan's lemma, one gets
\begin{align}\label{step2}
dt &=3\,t \,\phi_{1,1} +T_{{-1}}\, \omega+T_{{0}}\,\theta+T_{{1}}\,\eta,   \\
\del_{0,0} &=\frac{1}{3}\,T_{{1}}\,\omega+A_{{2}}\,\theta, \n \\
\del_{1,2} &=-3\,A_{{2}}\,\omega+2\,A_{{3}}\,\theta,  \n \\
\del_{0,1} &=-t\,\phi_{1,2} - \left( T_{{0}}+3\,t A_{{2}} \right) \,\omega
+2\,A_{{4}} \, \theta-3\,A_{{2}}\, \eta, \n
\end{align}
for coefficients $\{ \, T_{-1}, \, T_0, \, T_1 \}, \, \{ A_2, A_3, A_4 \}$.

\texttt{Step 3.} Differentiating $\, \del_{0,0} = \del_{2,2}$,
and after applying Cartan's lemma, one gets
\beq\label{step3}
\del_{0,2} = \left( t A_{{3}}+A_{{4}} \right) \omega+A_{{5}}\theta+A_{{3}}\eta, \n
\eeq
for a coefficient $\, A_5$.

\texttt{Step 4.} Differentiating the first equation of  \eqref{step2},
one gets
\begin{align}
d T_0 &=T_{{0}} \left( 4\,\phi_{{1,1}}+2\,\phi_{{2,2}} \right)
+T_{{-1}}\phi_{{1,2}}-T_{{1}}\phi_{{0,1}}+T_{{0,-1}}\omega+T_{{0,0}}\theta+T_{{1,0}}\eta,  \n \\
d T_1 &= T_{{1}} \, \left( 2\,\phi_{{1,1}}+\phi_{{2,2}} \right)
+3\,t  \, \phi_{{1,2}}+
T_{{1,-1}} \, \omega+T_{{1,0}} \, \theta+T_{1,1} \,\eta. \n
\end{align}
Differentiating  $\, \del_{0,0}  =\frac{1}{3}\,T_{{1}}\,\omega+A_{{2}}\,\theta$
with these equations, one gets $\,  T_{1,1}=12  \,A_2$, and
\beq
d A_2 = A_{{2}} \,  \left( \phi_{{1,1}}+2\,\phi_{{2,2}} \right)
+\frac{1}{3} \,T_{{1}} \,  \phi_{{1,2}}+ \left( \frac{1}{3} \,T_{{1,0}}+A_{{4}}-t A_{{3}} \right)  \, \omega
+A_{{2,0}} \,  \theta-A_{{3}} \, \eta. \n
\eeq

\texttt{Step 5.} Successive differentiation  of the rest of the equations in \eqref{step2}
implies that
\begin{align}
d A_3&=3\,A_{{3}}\phi_{{2,2}}-3\,A_{{2}}\phi_{{1,2}}+ \left( -\frac{3}{2}\,A_{{2,0}}-
\frac{1}{2}\,A_{{5}}+\frac{9}{2}\,{A_{{2}}}^{2}+A_{{3}}T_{{1}} \right) \omega+A_{{3,0}}\theta, \n \\
d A_4&=3\,A_{{4}} \left( \phi_{{1,1}}+\phi_{{2,2}} \right) +3\,A_{{2}}\phi_{{0,1}}
- \left( T_{{0}}+3\,t A_{{2}} \right) \phi_{{1,2}} \n \\
&\; \; - \left( T_{{1}}A_{{4}}
+3\,A_{{2}}T_{{0}}+\frac{9}{2}\,t {A_{{2}}}^{2}+\frac{3}{2}\,t A_{{2,0}}
+\frac{1}{2}\,T_{{0,0}}-\frac{1}{2}\,t A_{{5}} \right) \omega+A_{{4,0}}\theta+ \left(\frac{1}{2}\,A_{{5}}
-\frac{9}{2}\,{A_{{2}}}^{2}-\frac{3}{2}\,A_{{2,0}} \right) \eta, \n \\
d A_5&=A_{{5}} \left( 2\,\phi_{{1,1}}+4\,\phi_{{2,2}} \right) -3\,A_{{3}}\phi
_{{0,1}}+ \left( 3\,t A_{{3}}+3\,A_{{4}} \right) \phi_{{1,2}} \n \\
&\; \; - \left( -3\,A_{{3}}T_{{0}}-6\,A_{{2}} t A_{{3}}+6\,A_{{2}}A_{{4}}-A_{{4,0}}
-t A_{{3,0}} \right) \omega+A_{{5,0}}\theta+ \left( 6\,A_{{2}}A_{{3}}+A_{{3,0}} \right) \eta. \n
\end{align}
The fundamental structure equation \eqref{fundast} for $\, \del$
is now an identity.

\begin{theorem}\label{P2}
There exists a unique flat path geometry on $\, \PP^2$.
\end{theorem}
\Pf
We show that the deformation $\, \del$ vanishes.

Let $\, Z = (Z_0, \, Z_1, \, Z_2)$ be the  projective frame
of $\, SL_3 \C$ that satisfies the structure equation
$\,
d Z = Z\, \phi.
$
Let $\, Z^* = ( Z^0, \, Z^1, \, Z^2)^{t}$ be the dual frame
that  satisfies the structure equation
\beq\label{coframe}
d Z^* = - \phi  \, Z^*.\n
\eeq

Recall $\, \PP^2 = \PP(V), \, V = \C^3$.
Consider  a $\, Sym^3V^*$-valued function
\beq\label{sigma}
\sigma = t\, (Z^{{1}})^{3}+T_{{1}} \, (Z^{{1}})^{2} Z^{{2}}+
6\,A_{{2}} \,Z^{{1}}(Z^{{2}})^{2}-2\,A_{{3}} \, (Z^{{2}})^{3}.
\eeq
The structure equation    shows that
\beq\label{dsigma}
d \sigma \equiv 0, \mod \; \, \omega, \, \theta.
\eeq
$\, \sigma$ defines a map
\beq
\sigma: \, \PP^2 \to Sym^3V^*. \n
\eeq
Since $\, Sym^3V^*$ is a vector space, and $\, \PP^2$ is compact,
$\, \sigma$ is constant.

Equation \eqref{dsigma} contains the following terms.
\begin{align}
d \sigma &\equiv  (-3 t \, \omega - T_1 \, \theta) \, (Z^1)^2 Z^0
                        +( T_{-1} \, \omega + T_0 \, \theta) \,  (Z^1)^3,
                         \mod \; \, Z^2, \n \\
               &= 0.  \n
\end{align}
Note that the coefficients of $\, \del$
are polynomials in the   derivatives of $\, t$.
$\sq$

Let $\, \mathcal{E}$ be the rank two co-normal  bundle over $\, \PP^2 = \PP(V)$,
\beq
\mathcal{E} = \{ \, (\, [v], \, w) \in  \PP^2 \times \,  V^* \, | \; \, w(v) = 0  \, \}. \n
\eeq
$\, \mathcal{E}$ is generated by sections of $\, \{ \, Z^1, \, Z^2 \, \}$,
and $\, \sigma \in H^0(\PP^2,  Sym^3 \mathcal{E})$.
Let $\, T^*\PP^2$ be the cotangent bundle of $\, \PP^2$.
\begin{lemma}
\beq
\mathcal{E}  =T^*\PP^2  \otimes \mathcal{O}(1). \n
\eeq
\end{lemma}
\Pf
Let  $\, \mu = f_1 Z^1 + f_2 Z^2$ be a section of $\, \mathcal{E}$.
From the structure equation,
\begin{align}
df_1&\equiv   f_{{1}}\phi_{{1,1}}+ f_{{2}}\eta,  \n \\
df_2 &\equiv  f_{{2}}\phi_{{2,2}}+ f_{{1}}\phi_{{1,2}}, \mod \; \omega, \, \theta. \n
\end{align}
On the other hand,
let   $\, \alpha = F_1 \omega + F_2 \theta$ be a section of $\, T^*\PP^2$.
From the structure equation,
\begin{align}
dF_1&\equiv   F_{{1}}(2 \phi_{{1,1}}+\phi_{2,2}) + F_{{2}}\eta,  \n \\
dF_2 &\equiv   F_{{2}}( \phi_{{1,1}}+2 \phi_{2,2}) +  F_{{1}}\phi_{{1,2}}, \mod \; \omega, \, \theta. \n
\end{align}
$\,  \alpha$ corresponds to the $\, V^* \otimes V$-valued function
\beq
(F_1 \, Z^1 + F_2 \, Z^2) \otimes Z_0. \; \;  \sq  \n
\eeq

The analysis in this section is local.
Since $\, \sigma$, \eqref{sigma}, is a symmetric cubic differential
up to scaling by the elements of the line bundle $\, \mathcal{O}(3)$,
it defines a possibly degenerate 3-web.
\begin{proposition}\label{cubic}
Let $\, M \subset \PP^2$ be an open subset with the induced flat path geometry.
Suppose $\, M$ is given another  flat path geometry structure.
There exists a section $\, \sigma \in  H^0(M, \, Sym^3(T^*M) \otimes \mathcal{O}(3))$
such that the two path geometries coincide along
a generalized 3-web $\, \W_{\sigma} = \sigma^{-1}(0) \subset \PP(TM)$
defined by the zero locus of $\, \sigma$.
Two path geometries are isomorphic when $\, \sigma \equiv 0$.
\end{proposition}
When $\, \sigma \not \equiv 0$, and $\, \sigma$ has no repeated roots,
$\, \W_{\sigma} \subset \PP(TM)$ consists of three sections
$\, M \hook \PP(TM)$ defined by the base locus of $\, \sigma$.
It should be noted that
this does not imply   the associated 3-web of foliations  on $\, M$ is linear;
$\, \W_{\sigma} \subset \PP(TM)$  is not necessarily tangent
to the standard paths of $\, \PP(TM)$.

The 3-web $\, \W_{\sigma}$ is nonetheless not arbitrary.
An examination of the structure equation  shows that
the nontrivial part is $\, \{ \,A_{4,-1}, \, A_{5,-1} \, \}$,
the $\, \omega$-derivatives of $\, \{ \, A_4, \, A_5 \, \}$.
This observation  is the basis of
differential analysis toward Gronwall conjecture.

%-------------------------------------------------
%-------------------------------------------------
\section{Web curvature}\label{sec3}
Let $\, \W$ be  the   3-web  on a surface $\, M$
defined by three 1-forms
$\, \omega^i, \, i=1, \, 2, \, 3$, such that $\, \omega^i \w \omega^j \ne 0$ for $\, i \ne j$.
Up to scaling, one may arrange so that
\beq\label{sum0}
\omega^1 + \omega^2 + \omega^3 = 0.
\eeq
The principal $\, GL_2\C$-frame bundle of $\, M$ can be reduced to
a $\, G \subset GL_2\C$-bundle F$\to M$  by \eqref{sum0}, where
\beq
G = \{ \, \lambda  \begin{pmatrix}
1 & 0 \\
0 & 1 \\
\end{pmatrix} \;   | \; \; \lambda \ne 0 \, \}. \n
\eeq
We continue to use  $\, \omega^i$
to denote the tautological semi-basic 1-forms on F.

Let $\, \rho$ be the unique connection 1-form on F   that satisfies
\beq
d \omega^i + \rho \w \omega^i = 0, \; \; i=1, \, 2, \, 3.  \n
\eeq
\emph{Web curvature} $\, K$ of $\, \W$ is the relative invariant of
the induced  $\, G$-structure F defined by
\beq\label{curvature}
d \rho = K \, \omega^1 \w  \omega^2.
\eeq
The functional relations among $\, K$  and its successive derivatives
are the  basic  local invariants of $\, \W$.

Let $\, \W$ be the 3-web on an open subset  $\, M \subset \PP^2$
which arises as the intersection of two distinct  flat path geometries.
Let $\, \sigma   \in   H^0(M,  \, Sym^3(T^*M) \otimes \mathcal{O}(3))$
be the associated section.
The condition that $\, \W$ is linear is expressed as a conformally invariant
first order differential equation for $\, \sigma$, Lemma \ref{sigmalinear}.
Under this linearity condition,
we derive an explicit formula  for the web curvature, \eqref{webcurvature},
which depends on the fifth order jet of the deformation function $\, t$.

Furthermore, successive differentiation via the over-determined PDE machinery
allows one to close up the structure equation for $\, \W$  at order eight.
In the course of computation,
a set of structure coefficients will be normalized
by the more or less standard frame adaptation.
This is relevant for our purpose
where it is necessary to keep the  computations under  manageable size.

As in Section \ref{sec2},
the notational convention for the  covariant derivatives of a function is;
\beq
d f \equiv f_{-1} \omega + f_0 \, \theta + f_1\, \eta, \mod \; \phi, \n
\eeq
where  (mod $ \; \phi)  $ means
 (mod $  \; \phi_{1,1}, \, \phi_{2,2}, \, \phi_{0,1}, \, \phi_{1,2}, \, \phi_{0,2}) $.

\subsection{Linearity}\label{sec31}
In this subsection, we determine the compatibility equation
for the intersection of two distinct  flat path geometries
to be a linear 3-web.

Let $\, \W$ be a linear 3-web on an open subset $\, M \subset \PP^2$
with respect to the standard flat path geometry $\, \Lambda$ of $\, \PP^2$.
Let  $\, \Lambda'$  be  another distinct flat path geometry structure  on $\, M$.
Let $\, \sigma \in   H^0(M, \, Sym^3(\mathcal{E}))$ be the associated  section.
Assume $\, \Lambda'$ also  linearizes $\, \W$, or equivalently
each leaf of the foliations of $\, \W$ is a part of a path of $\, \Lambda'$.
It is clear that $\, \W$ must   coincide with $\, \W_{\sigma}$.
Since $\, \W= \W_{\sigma}$ is linear,
this puts further differential geometric constraints on $\, \sigma$,
see remark below Proposition \ref{cubic}.
\begin{remark}
This observation implies   that a linear  $\, d$-web for  $\, d \geq 4$
has a unique linearization up to projective transformation.
\end{remark}

Let $\, \sharp : \mathcal{E}   \to T^*M$ be an isomorphism defined by
a nonzero section of $\, \mathcal{O}(-1)$.
Let $\, \flat: Sym^3(\mathcal{E}) \otimes T^*M \to Sym^4(\mathcal{E}) $
denote the associated  symmetrization.
\begin{lemma}\label{sigmalinear}
3-web $\, \W_{\sigma} = \sigma^{-1}(0)$ is linear when
\beq\label{slinear}
(d \sigma)^{\flat} \equiv 0, \mod \; \, \sigma.
\eeq
\end{lemma}
\Pf
Note first that this equation is well defined independent of the section of
$\, \mathcal{O}(-1)$ defining the isomorphism  $\, \sharp$.
A different section only scales $\, \sharp$, and does not affect  \eqref{slinear}.

Decompose $\, \sigma = \sigma_1 \, \sigma_2 \, \sigma_3$,
$\, \sigma_i \in H^0(M, \, \mathcal{E})$.
Let $\, \alpha_i = (\sigma_i)^{\sharp}$ be the corresponding 1-form.
By definition, each $\, \sigma_i$-curve is linear when
\beq
d \sigma_i \equiv 0, \mod \; \, \sigma_i; \, \alpha_i.
\eeq
It can be checked by elementary computation that this is equivalent to
$\, (d \sigma)^{\flat}$ being divisible by $\, \sigma_1, \, \sigma_2$, and  $\, \sigma_3$.
$\sq$

Set
\begin{align}
(Z^1)^{\sharp} &= \omega, \n \\
(Z^2)^{\sharp} &= \theta. \n
\end{align}
From \eqref{coframe}  and the structure equation,
\begin{align}
(d \sigma)^{\flat} + 3 \, Z^0 \, \sigma
&=T_{{-1}}\, (Z^{{1}})^{4}+ \left( T_{{1,-1}}+T_{{0}} \right) (Z^{{1}})^{3} Z^{{2}}
+ \left( 3\,T_{{1,0}}+6\,A_{{4}}-6\,t A_{{3}} \right) (Z^{{1}})^{2}(Z^{{2}})^{2} \n \\
&\; + \left( A_{{5}}-9\,{A_{{2}}}^{2}-2\,T_{{1}}A_{{3}}+9\,A_{{2,0}} \right)
Z^{{1}} (Z^{{2}})^{3}-2\,A_{{3,0}} \, (Z^{{2}})^{4}.  \n
\end{align}
\eqref{slinear} is equivalent to
\begin{align}\label{3eq}
-tA_{{3}}T_{{1,-1}}+{t}^{2}A_{{3,0}}&=
-T_{{1}}T_{{-1}}A_{{3}}+tA_{{3}}T_{{0}},   \\
-3\,t A_{{3}}T_{{1,0}} + t \, T_{{1}}A_{{3,0}}
&=-6\,{t}^{2}{A_{{3}}}^{2}+6\,tA_{{3}}A_{{4}}-6\,A_{{2}}T_{{-1}}A_{{3}}, \n \\
-9\,t \,A_{{3}}A_{{2,0}} + 6 \, t\,A_{{2}}A_{{3,0}}&=
-9\,tA_{{3}}{A_{{2}}}^{2}+2\,{A_{{3}}}^{2}T_{{-1}}
+tA_{{3}}A_{{5}}-2\,t{A_{{3}}}^{2}T_{{1}}. \n
\end{align}

Solving this system of equations for $\, \{ \, T_{1,0}, \, A_{2,0}, \, A_{3,0} \}$,
\begin{align}\label{eqq1}
T_{1,0} &=  \,{\frac {6\,T_{{-1}}t A_{{2}}-6\,{t}^{2}A_{{4}}+6\,{t}^{3}A_{{3}}
-{T_{{1}}}^{2}T_{{-1}}+t  T_{{0}}T_{{1}}+t T_{{1,-1}}T_{{1}}}{3{t}^{2}}},     \\
A_{2,0} &=  \,{\frac {-2\,A_{{3}}T_{{-1}} t +2\,T_{{1}}{t  }^{2}A_{{3}}
-6\,A_{{2}}T_{{1}}T_{{-1}}-{t}^{2}A_{{5}}+9\,{t}^{2}{A_{{2}}}^{2}
+6\,t A_{{2}}T_{{1,-1}}+6\,t A_{{2}}T_{{0}}}{9{t}^{2}}},  \n \\
A_{3,0} &= {\frac {A_{{3}} \left( -T_{{1}}T_{{-1}}+t T_{{0}}+t T_{{1,-1}} \right) }{{t}^{2}}}.  \n
\end{align}
For a later purpose, solving the system of equations \eqref{3eq}
for $\, \{ \, T_{-1}, \, A_{2,0}, \, A_{3,0} \}$,
\begin{align}\label{eqq2}
T_{-1} &= -{\frac {t \left( 6\,{t}^{2}A_{{3}}-6\,A_{{4}} t
+T_{{1}}T_{{1,-1}}-3\,T_{{1,0}} t +T_{{1}}T_{{0}} \right) }{6\,t A_{{2}}-{T_{{1}}}^{2}}},\\
A_{2,0} &=  \,{\frac{1}{9 ({6\,t A_{{2}}-{T_{{1}}}^{2}} )} }
 (12\,{t}^{2}{A_{{3}}}^{2}+54\,t {A_{{2}}}^{3}
-12\,t A_{{3}}A_{{4}}+48\,t A_{{2}}T_{{1}}A_{{3}}-6\,t A_{{3}}T_{{1,0}}
-36\,A_{{2}}A_{{4}}T_{{1}} \n \\
&-6\,t A_{{2}}A_{{5}}-9\,{A_{{2}}}^{2}{T_{{1}}}^{2}
+2\,T_{{1}}
A_{{3}}T_{{0}}-2\,A_{{3}}{T_{{1}}}^{3}+36\,{A_{{2}}}^{2}T_{{1,-1}}
+2\,T_{{1}}A_{{3}}T_{{1,-1}}+36\,T_{{0}}{A_{{2}}}^{2} \n \\
&-18\,T_{{1}}T_{{1,0}}
A_{{2}}+A_{{5}}{T_{{1}}}^{2}),   \n \\
A_{3,0} &=3\,{\frac {A_{{3}} \left( 2\,T_{{1}} t A_{{3}}
-2\,T_{{1}}A_{{4}}-T_{{1,0}}T_{{1}}+2\,A_{{2}}T_{{1,-1}}+2\,A_{{2}}T_{{0}} \right) }
{6\,t A_{{2}}-{T_{{1}}}^{2}}}.  \n
\end{align}
\begin{remark}\label{310}
It will be shown in   Section \ref{sec32} that both
 $\, T_1$, and $\, 9\, A_2^2 + 2\,T_1 A_3$ are nonzero  on $\, t^{-1}(0)$.
Likewise,  some of the denominators in the expressions
are invertible elements in the appropriate local ring,
and the divisions  make sense.
 \end{remark}

The   web curvature of $\, \W_{\sigma}$
will be computed  by restricting to the zero locus of the deformation function $\, t$.
Let us introduce the standard notation for local ring $\,  \mathcal{O}(t)$, \cite[p10]{GH}.
Some of the expressions  below for example
have too many terms to be written down completely.
These will be written   in mod $\, \mathcal{O}(t)$.

 Set
\beq\label{T-11}
 d T_{1,-1} = T_{{1,-1}} \left( 4\,\phi_{{1,1}}+2\,\phi_{{2,2}} \right) -3\,t\phi_{{0
,2}}-2\,T_{{1}}\phi_{{0,1}}+3\,T_{{-1}}\phi_{{1,2}}+T_{{1,-1,-1}}
\omega+T_{{1,-1,0}}\theta+T_{{1,-1,1}}\eta.
\eeq
Exterior derivatives $\, d(d(T_1)), \, d(d(A_2))$  with \eqref{eqq1},
and then evaluating with \eqref{eqq2} for $\, T_{-1}$
give
\beq
T_{1,-1,1}  =5\,T_{{1,0}}+12\,A_{{4}}-12\,tA_{{3}},  \n
\eeq
and
\begin{align}\label{eqq3}
T_{-1,-1}  &\equiv  \frac{1}{3 T_{{1}}^{7}} ( {3\,{T_{{1}}}^{6}T_{{1,-1,-1}}
+6\,{T_{{0}}}^{2}{T_{{1}}}^{5}+6\,T_{{0}}T_{{1,-1}}{T_{{1}}}^{5}}) t, \mod \; \mathcal{O}(t^2),   \\
-T_{{1}}T_{{0,-1}}+T_{{0}}T_{{1,-1}}+{T_{{0}}}^{2} &\equiv
-  \frac{1}{3 T_1^7} (  36\,A_{{2}}{T_{{1}}}^{5}T_{{1,-1}}T_{{0}}
+3\,{T_{{1}}}^{7}T_{{0,0}}+36\,A_{{2}}{T_{{1}}}^{5}{T_{{0}}}^{2}
-21\,{T_{{1}}}^{6}T_{{1,0}}T_{{0}}  \n \\
&\; \; \; \; +3\,{T_{{1}}}^{7}T_{{1,-1,0}}
-3\,{T_{{1}}}^{6}T_{{1,-1}}T_{{1,0}}-36\,{T_{{1}}}^{6}A_{{4}}T_{{0}}   )  t,
\mod \; \mathcal{O}(t^2). \n
\end{align}

Taking exterior derivative  $\,  d(d(A_3))$ with these relations, and solving for $\, A_{5,0}$,
one gets a compatibility equation for $\, \W_{\sigma}$ to be linear.
\begin{align}\label{comp1}
A_{5,0} &\equiv
\frac{2}{T_1^4}\, ( -18\,A_{{2}}A_{{4}}{T_{{1}}}^{2}T_{{1,0}}+36\,{A_{{2}}}^{2}
T_{{1}}T_{{1,-1}}T_{{1,0}}-18\,{A_{{2}}}^{2}T_{{1}}T_{{1,0}}T_{{0}}+2
\,A_{{3}}{T_{{1}}}^{2}T_{{1,0}}T_{{0}}+8\,A_{{3}}{T_{{1}}}^{2}T_{{1,-1
}}T_{{1,0}} \\
&+108\,{A_{{2}}}^{3}{T_{{0}}}^{2}
-108\,{A_{{2}}}^{3}{T_{{1,-1}}}^{2}+2\,{T_{{1}}}^{3}A_{{5}}T_{{1,0}}
+4\,{T_{{1}}}^{3}A_{{5}}A_{{4}}-36\,A_{{2}}{A_{{4}}}^{2}{T_{{1}}}^{2}
-12\,A_{{3}}A_{{2}}T_{{0}}T_{{1,-1}}T_{{1}} \n \\
&-18\,{A_{{2}}}^{2}T_{{1,-1,0}}{T_{{1}}}^{2}
+7\,A_{{3}}A_{{4}}{T_{{1}}}^{4}+36\,{A_{{2}}}^{2}T_{{0,0}}{T_{{1}}}^{2}+12\,A_{{3}}A
_{{2}}{T_{{0}}}^{2}T_{{1}}-24\,A_{{3}}A_{{2}}{T_{{1,-1}}}^{2}T_{{1}}+
30\,A_{{3}}{T_{{1}}}^{3}A_{{2}}T_{{0}} \n \\
&+36\,T_{{1}}{A_{{2}}}^{2}A_{{4}}
T_{{0}}+144\,T_{{1}}{A_{{2}}}^{2}A_{{4}}T_{{1,-1}}-6\,T_{{1,-1}}A_{{5}
}{T_{{1}}}^{2}A_{{2}}-6\,T_{{0}}A_{{5}}{T_{{1}}}^{2}A_{{2}}+12\,A_{{3}
}T_{{1,-1,-1}}{T_{{1}}}^{2}A_{{2}} \n \\
&+8\,A_{{3}}{T_{{1}}}^{2}T_{{0}}A_{{4
}}+20\,A_{{3}}{T_{{1}}}^{2}T_{{1,-1}}A_{{4}}+5\,A_{{3}}{T_{{1}}}^{3}T_
{{0,0}}-4\,A_{{3}}{T_{{1}}}^{3}T_{{1,-1,0}}+6\,{T_{{1}}}^{3}A_{{2}}A_{
{4,0}}+54\,T_{{1}}{A_{{2}}}^{3}T_{{1,-1,-1}} \n \\
&+162\,{T_{{1}}}^{2}{A_{{2}
}}^{3}T_{{0}}+36\,{T_{{1}}}^{3}{A_{{2}}}^{2}A_{{4}} ),
\mod \; \mathcal{O}(t). \n
\end{align}

The computation of  exterior derivatives $\,  d(d(T_1))$, $\,  d(d(A_2))$, and $\,  d(d(A_3))$
using \eqref{eqq1}
are performed away from the zero locus of $\, t$.
The compatibility equation  thus obtained admits an
analytic continuation  to $\, t^{-1}(0)$
to yield  \eqref{comp1}.

The idea of this argument   is the elementary residue theorem.
It can be described in the current setting as follows.
Let $\, f$ be an analytic function on $\, \PP(TM) \to M$ which vanishes on   $\,  t^{-1}(0)$.
Since $\, T_1$,  the derivative of defining function $\, t$
in the fiber direction($\, \eta$-derivative),
is nonzero at $\, t^{-1}(0)$, Remark \ref{310},
$\, f$ is divisible by $\, t$ and one may write
$\,
f = f^{(1)} \, t + f^{(2)} \, t^2 + \, ... \, .
$
Taking $\, \eta$-derivative, one gets
\beq
f_1 = \frac{\partial f}{\partial \eta} \equiv  f^{(1)} \, T_1, \mod \; \mathcal{O}(t). \n
\eeq
The coefficient $\, f^{(1)} \equiv  \frac{f_1}{T_1} \mod \; \mathcal{O}(t)$ is obtained
by differentiating $\, f$ once in the fiber direction.

In \eqref{eqq3},
one could solve for $\, T_{1,-1,0}, \, T_{1,-1,-1}$ instead, away from $\, t^{-1}(0)$.
Together with \eqref{comp1} for $\, A_{5,0}$, this implies that
the set of remaining independent variables whose covariant derivatives are not determined
at this stage is
$\, \{\, A_{4,0}, \, T_{0,0}, \, T_{0,-1}, \, T_{-1,-1} \}$.

Successive differentiation  suggests that
the structure equations for the linear 3-web $\, \W_{\sigma}$
eventually close up.
Moreover, there appear to exist a triality of flat path geometries.
Since this is not  directly related to Gronwall conjecture in our treatment,
let us give a sketch of ideas.

Let $\, \langle \, \theta \, \rangle^{\perp} \subset T(\PP(TM))$ be the canonical contact 2-plane field.
Three corank one Pfaffian systems
$\,  \langle \, \theta, \, \omega\, \rangle, \,  \langle \, \theta, \, \eta \,\rangle,
 \,  \langle \, \theta, \, \eta+ t\, \omega \, \rangle$
define a set of three transversal Legendrian foliations
$\, \mathcal{F}^0, \, \mathcal{F}^1, \, \mathcal{F}^{-1}$ respectively away from $\, t^{-1}(0)$.
The assumption is that the two pairs
$\, ( \mathcal{F}^0, \, \mathcal{F}^1)$, and
$\, ( \mathcal{F}^0, \, \mathcal{F}^{-1})$ are projectively flat.
An analysis of prolonged structure equations seems to indicate that
in this case
the third pair of  Legendrian foliations
$\, ( \mathcal{F}^1, \, \mathcal{F}^{-1})$ is also projectively flat.

%%%%%%%%%%%%%%%%%%%%%%%%%%%%%%%%%%%%%%
%%%%%%%%%%%%%%%%%%%%%%%%%%%%%%%%%%%%%%
\subsection{Web curvature}\label{sec32}
For the rest of the paper,
we restrict our attention to the zero locus $\, t^{-1}(0)$.
In this subsection,
we derive a formula for the web curvature of  $\, \W_{\sigma}$
from the structure equation obtained in Section \ref{sec31}.

Let $\, L_0 \subset \PP(TM)$ be the fiber of the $\, \PP^1$-bundle
$\,  \PP(TM) \to M$ at  a reference point x$_0 \in M \subset \PP^2$.
Assuming locally,
3-web $\, \W_{\sigma}$
can be identified with the disjoint union of three surfaces
$\, M^a \subset \PP(TM), \, a=1, \, 2, \, 3$,
which is the zero locus of the deformation function $\, t$, \eqref{sigma}.
Equation \eqref{step2} implies that
$\, t$ represents a section of $\, \mathcal{O}(3)$
when restricted to a $\, \PP^1$-fiber of  $\,  \PP(TM) \to M$
(we omit the details).
Hence $\, T_1$,
the $\, \eta$-derivative, or the derivative of $\, t$  in the fiber direction,
does not vanish at  $\, M^a$,  and
$\, L_0$ intersects each $\, M^a$  transversally.
This  implies in particular  that any section of
$\, \{ \, \omega, \, \theta \, \}$ is a coframe on $\, M^a$.

Consider the incidence double fibration.

\begin{picture}(300,66)(-67,-5)\label{11double}
\put(147,40){$\PP(TM)$}
\put(175,25){$\searrow$}
\put(145,25){$\swarrow$}
\put(184,10){$M^* \subset (\PP^2)^*$ }
\put(105,10){ $\PP^2 \supset M$ }
\put(187,30){$\pi_2$}
\put(133,30){$\pi_1$}
\end{picture}

\noi
Under the projection $\, \pi_2$,
each $\, M^a$ is mapped to an immersed curve $\, \Gamma^a \subset  M^*$.
The image $\, \pi_2(\pi_1^{-1}(\mbox{x})) = L_{\mbox{x}}$
is the dual line of    $\mbox{x} \in M$,
which intersects $\, \cup \, \Gamma^a$ transversally
at three points.

The paths of the standard path geometry on $\, M \subset \PP^2$
are integral curves of the Pfaffian system $\, \langle \, \theta, \, \eta \, \rangle$, Theorem \ref{defining}.
Since $\, \W_{\sigma}$ is linear,
$\, \eta \w \theta = 0$ on $\, M^a$.
The structure equation \eqref{step2} for $\, t$
restricted to $\, t^{-1}(0)$  becomes
\begin{align}\label{T1}
T_{-1}&=0,   \\
\eta      &= - \frac{T_0}{T_1}  \, \theta  , \quad  \mbox{on} \; \, M^a. \n
\end{align}
Differentiating this equation, one gets
\begin{align}\label{T01}
T_{-1,-1}&=0,   \\
T_{0,-1}&=(T_{1,-1}+T_0) \, \frac{T_0}{T_1}, \; \; \; \mbox{on} \; \, M^a.   \n
\end{align}
This agrees with \eqref{eqq3}.

For definiteness, fix  a single section  $\, M^1 \subset \PP(TM)$.
Since $\, \W = \W_{\sigma}$ is a 3-web,
$\, \sigma$ has three distinct roots on each fiber $\, \PP(TM) \to M$.
By definition,
\beq
\sigma |_{M^1} = ( T_{{1}} \, (Z^{{1}})^{2}
+ 6\,A_{{2}} \,Z^{{1}} Z^{{2}} -2\,A_{{3}} \, (Z^{{2}})^{2})\,  Z^{{2}}, \n
\eeq
and we must have
\beq\label{9A2}
 9 \, A_2^2 + 2\, T_1  A_3 \ne 0.
\eeq

%%%%%%%%%%%%%%%%%%%%%%%%%%%%%%%%%%%%%%
%%%%%%%%%%%%%%%%%%%%%%%%%%%%%%%%%%%%%%
\subsubsection{Normalization}\label{sec321}
From the structure equations of Section \ref{sec2},
\begin{align}
d T_0 &\equiv T_{{0}} \left( 4\,\phi_{{1,1}}+2\,\phi_{{2,2}} \right)
                        -T_{{1}}\phi_{{0,1}}, \mod \;   \, \omega, \, \theta, \, \eta, \n \\
d A_2 &\equiv  A_{{2}} \,  \left( \phi_{{1,1}}+2\,\phi_{{2,2}} \right)
                          +\frac{1}{3} \,T_{{1}} \,  \phi_{{1,2}}, \mod \;  \, \omega, \, \theta, \, \eta, \n \\
d T_{1,0} &\equiv  T_{{1,0}} \left( 3\,\phi_{{1,1}}+ 3 \,\phi_{{2,2}} \right)
                        -T_{{1}}\phi_{{0,2}}, \mod \; \, A_2; \, \phi_{1,2},   \, \omega, \, \theta, \, \eta, \n \\
d ( \frac{2 A_3}{T_1} )  &\equiv   \frac{2 A_3}{T_1} \left( -2 \,\phi_{{1,1}}+2 \,\phi_{{2,2}} \right),
                      \mod \;   \, A_2; \, \omega, \, \theta, \, \eta. \n
\end{align}
By \eqref{9A2},
one may use the group action that  corresponds to
$\, \{ \, \phi_{0,1}, \, \phi_{1,2}, \, \phi_{0,2}, \, (\phi_{1,1}-\phi_{2,2}) \, \}$
to  normalize
\begin{align}
T_0  &=  0, \n \\
A_2 &=   0, \n \\
T_{1,0} &=   0, \n \\
\frac{2 A_3}{T_1} &=  1. \n
\end{align}

Derivatives of these relations imply  that
\begin{align}
\phi_{0,1}                  &={\frac {T_{{0,0}}}{T_{{1}}}}\, \theta, \n \\
\phi_{1,2}                  &=- \,{\frac {3 A_{{4}}}{T_{{1}}}}\, \omega
+ {\frac {  ( T_{{1,-1}} - {T_{{1}}}^{2} +A_{{5}} )}{3 T_{{1}}}}\, \theta, \n \\
\phi_{0,2}                  &=- \,{\frac {  ( -8\,T_{{1,-1}}A_{{4}}-5\,T_{{0,0}}T_{{1}}+A_{{5,0
}}T_{{1}}-8\,A_{{5}}A_{{4}}-7\,A_{{4}}{T_{{1}}}^{2}  ) }{{4 T_{{1}}}^{2}}}\, \omega \n \\
&+ \,{\frac {  ( -30\,A_{{4,0}}T_{{1}}-3\,T_{{1,-1,-1}}T_{{1}}
+2\,{T_{{1}}}^{4}-2\,{T_{{1}}}^{2}A_{{5}}+5\,{T_{{1,-1}}}^{2}
+108\,{A_{{4}}}^{2}+5\,T_{{1,-1}}A_{{5}}-{T_{{1}}}^{2}T_{{1,-1}} ) }
{{9 T_{{1}}}^{2}}}\, \theta, \n \\
\phi_{1,1}-\phi_{2,2}&= \,{\frac {  ( -T_{{1,-1}}+{T_{{1}}}^{2}-A_{{5}}  ) }
{3 T_{{1}}}}\, \omega + {\frac {3 A_{{4}}}{T_{{1}}}} \, \theta,    \n
\end{align}
and that
\begin{align}\label{T00}
dT_{0,0} &= 9\,T_{{0,0}}\phi_{{1,1}} - \,{\frac {T_{{0,0}}
( -7\,T_{{1,-1}}+4\,{T_{{1}}}^{2}
-4\,A_{{5}}  )}{3\,T_{{1}}}} \omega+{\frac {  ( -12\,T_{{0,0}}A_{{4}}
+B_{{1}}T_{{1}}  )}{T_{{1}}}} \theta,   \\
T_{1,-1,0} &=- \,{\frac {-5 \,T_{{0,0}}T_{{1}}+A_{{5,0}}T_{{1}}-8\,A_{{5}}A_{{4}}-
20\,T_{{1,-1}}A_{{4}}-7\,A_{{4}}{T_{{1}}}^{2}}{4\,T_{{1}}}}, \n
\end{align}
for a new variable  $\, B_1$.

Exterior derivatives $\, d(d(T_1))$, $\, d(d(A_2))$,  $\, d(d(A_3))$, $\, d(d(A_4))$, $\, d(d(A_5))$,
and $\, d(d(T_{1,-1}))$
with these relations imply
\begin{align}\label{T40}
dA_{4,0} &= 9 \,  A_{4,0} \phi_{1,1} + A_{4,0,-1} \, \omega + B_2 \, \theta,\,   \\
dA_{5,0} &=  9 \, A_{5,0} \phi_{1,1} + A_{5,0,-1} \, \omega + B_3 \, \theta,\, \n \\
dT_{1,-1,-1} &=  9 \, T_{1,-1,-1} \phi_{1,1} + T_{1,-1,-1,-1} \, \omega + T_{1,-1,-1,0} \, \theta, \n
\end{align}
for new variables $\, B_2, \, B_3$,
where the coefficients $\, A_{4,0,-1}, \, A_{5,0,-1}, \, T_{1,-1,-1,-1}, \, T_{1,-1,-1,0}$
are polynomials in the known variables $\, T_i, \, T_{i,j}, T_{i,j,k}, A_i, \, A_{i,j}, B_j$.
Due to their lengths, the exact expressions are postponed to Appendix.

Note that we use a slightly different formulation compared to Section \ref{sec31}.
Here $\, A_{5,0}$ is an independent coefficient,  check \eqref{comp1}.

\subsubsection{Web curvature}\label{sec322}
$\, \mathcal{E} \subset M \times V^*$ is by definition
a sub-bundle of a trivial bundle.
The section  $\, \sigma \in H^0(M, Sym^3(\mathcal{E}) )$
can  be considered  as  a $\, Sym^3  V^*$-valued function on $\, M$.
For the purpose of differential analysis,
$\, \sigma$ is freely pulled back to $\, M^1$.

From the normalization,   $\, \sigma$ decomposes on $\, M^1 \subset t^{-1}(0)$
as a product of monomials
\beq
\sigma = T_1\, (Z^1+Z^2) (Z^1-Z^2) \, Z^2.
\eeq
$\, d Z_0 \equiv Z_1 \, \omega + Z_2 \, \theta, \mod \; Z_0$,
and since $\, (Z_0, Z_1, Z_2)$ is dual to $\, (Z^0, Z^1, Z^2)$,
$\, \W_{\sigma}$ is the 3-web on $\, M^1$ defined by
the base locus of the symmetric cubic differential
\beq
\Psi =   (\omega+\theta) \circ (\omega-\theta) \circ \theta. \n
\eeq

Set
\begin{align}\label{sigmaformula}
\sigma_1 &= Z^1+Z^2,   \\
\sigma_2 &= Z^1-Z^2,  \n \\
\sigma_3 &=T_1 \, Z^2, \n
\end{align}
and
\begin{align}
\alpha_1 &= \omega+\theta,  \n \\
\alpha_2 &= \omega-\theta,  \n \\
\alpha_3 &=  \, \theta. \n
\end{align}

By assumption,  an integral curve of each $\, \sigma_i$, $\, i=1, \, 2$,
projects to a part of a line in $\, \PP^2$.
This   translates to the equation
\beq
d \sigma_i \equiv 0, \, \mod \; \sigma_i; \, \alpha_i. \n
\eeq
Direct computation shows that this is equivalent to
\begin{align}\label{linearcond}
A_{2,0} &=\frac{-T_{{1,-1}}+{T_{{1}}}^{2}-A_{{5}}}{9},   \\
A_{3,0} &=3\, A_4. \n
\end{align}
One may check that this set of equations agree with
\eqref{eqq2} restricted to $\, M^1$.

Let $\, \rho$ be the unique 1-form such that
\beq
d \alpha_i + \rho \w \alpha_i =0, \; \, i = 1, \, 2, \, 3.  \n
\eeq
Direct computation using \eqref{linearcond} shows that
\begin{align}
\rho &=3\,\phi_{{1,1}}- \,{\frac {2 \left( -T_{{1,-1}}+{T_{{1}}}^{2}-A_{{5}}
 \right)}{3\, T_{{1}}}} \omega. \n
\end{align}
Differentiating the connection form $\, \rho$,   the web curvature $\, K$ is given by
\begin{align}
d \rho  =  K \, \omega \w \theta, \n
 \end{align}
 where
 \begin{align}\label{webcurvature}
 K &= -\,{\frac {  ( 4\,T_{{1,-1}}A_{{4}}+A_{{5,0}}T_{{1}}+5\,A_{{4}}
{T_{{1}}}^{2}+3\,T_{{0,0}}T_{{1}}-8\,A_{{5}}A_{{4}} \ )
 }{{2 T_{{1}}}^{2}}}.
\end{align}

\subsection{Prolongation}
In this subsection, we present the prolongation steps of
how  the structure equations for a linear 3-web
admitting two distinct linearizations
close up at order eight.

Due to their lengths, the exact expressions for the coefficients
$\, \{ \,     B_{1,-1}$, $B_{2,-1}$, $B_{3,-1}$, $B_{3,0}$, $B_{4,-1}$,
$B_{5,-1}$, $B_{5,0}$, $B_{6,-1}$, $B_{6,0}   \, \}$ below
are postponed to Appendix.

Note the relations
\begin{align}
A_{4,-1}&\equiv -\frac{1}{2} \, B_1, \n \\
A_{5,-1}&\equiv                       \, B_2, \n \\
T_{1,-1,-1,-1}&\equiv \frac{1}{4}(3\, B_3+5 \, B_1), \n \\
T_{1,-1,-1,0} &\equiv -\frac{1}{4} \, B_2, \mod \; \, \mbox{lower order terms}. \n
\end{align}
Differentiating \eqref{T00}, \eqref{T40},
and applying Cartan's lemma, one gets
\begin{align}\label{B1}
d B_1 &=12 \, B_1 \, \phi_{1,1} + B_{1,-1}\,  \omega + B_4 \, \theta,   \\
d B_2 &=12 \, B_2 \, \phi_{1,1} +B_{2,-1}\,  \omega + B_5 \, \theta,  \n  \\
d B_3 &=12 \, B_3 \, \phi_{1,1} +B_{3,-1}\,  \omega + B_{3, 0} \, \theta,  \n
\end{align}
for  new coefficients $\, B_4, \, B_5$.
The coefficients $\, B_{1,-1}, \, B_{2,-1}, \, B_{3,-1}, \, B_{3,0}$
are polynomials in the known variables, including $\, B_4, \, B_5$.

Note the relations
\begin{align}
B_{2,-1}&\equiv -\frac{1}{2} \, B_4 - \frac{10}{3}\, T_1\, B_2, \n \\
B_{3,-1}&\equiv B_5 + \frac{1}{2}\, T_1  \, B_1 -  \frac{8}{3}\, T_1  \, B_3, \n \\
B_{3,0} &\equiv -\frac{3}{2}\,  B_4 - 5\, T_1  \, B_2, \mod \; \, \mbox{lower order terms}. \n
\end{align}
Differentiating \eqref{B1}, one gets
\begin{align}\label{B4}
d B_4 &=15 \, B_4 \, \phi_{1,1} +B_{4,-1}\,  \omega + B_6 \, \theta,     \\
d B_5 &=15 \, B_5 \, \phi_{1,1} +B_{5,-1}\,  \omega + B_{5, 0} \, \theta,  \n
\end{align}
for  a new coefficient  $\, B_6$.
The coefficients $\, B_{4,-1}, \,  B_{5,-1}, \, B_{5,0}$
are polynomials in the known variables, including $\, B_6$.

Note the relations
\begin{align}
B_{4,-1}&\equiv -\frac{8}{3} \, T_1 \, B_4, \n \\
B_{5,-1}&\equiv - \frac{1}{2}\,  B_6 - 4\, T_1  \, B_5, \n \\
B_{5,0} &\equiv T_1\,  B_4 + \frac{26}{9} \, T_1^2  \, B_2, \mod \; \, \mbox{lower order terms}. \n
\end{align}
Differentiating \eqref{B4}, one finally gets
\beq\label{B6}
d B_6  =18 \, B_6 \, \phi_{1,1} +B_{6,-1}\,  \omega + B_{6, 0} \, \theta,
\eeq
where   the coefficients $\,  B_{6,-1}, \, B_{6,0}$
are polynomials in the known variables.

The structure equation for the linear 3-web $\, \W_{\sigma}$ closes up at this step.
Exterior derivative of  $\, d B_6$  yields a universal integrability condition,
\beq\label{X9}
0=d ( d ( B_6) ) = Eq_6 \, \omega \w \theta.
\eeq
$\, Eq_6$ is a polynomial with large number of terms.
Vanishing of $\, Eq_6$, and its successive derivatives
are necessary compatibility conditions
for the 3-web $\, \W_{\sigma}$ to admit two distinct linearizations.
Since $\, Eq_6$ will not be used, the exact expression shall be omitted.

%-------------------------------------------------
%-------------------------------------------------
\section{Gronwall conjecture}\label{sec4}
Web curvature \eqref{webcurvature} is a fifth order differential invariant of
the deformation function $\, t$,
whereas the universal integrability equation $\, Eq_6$
depends on the eighth order jet of $\, t$.

One way to approach Gronwall conjecture would be
by using the polynomial compatibility equations
obtained by successive differentiation of  $\, Eq_6$.
By solving for the higher order terms, $\, B_i$,
in terms of the lower order terms, $\, T_i, \, T_{i,j}, \, A_i, \, A_{i,j}$,
one may gradually lower the order of the compatibility equations
so that they become equations among the lower order terms only.
The simultaneous vanishing of  these compatibility equations
could then imply, for instance,
that the possible values of
$\, T_i, \, T_{i,j}, \, A_i, \, A_{i,j}$ are at most finite.
Hence their derivatives must be zero,
which would force the web curvature to vanish.

This method requires either solving a sequence of polynomial equations,
or showing that they have only finitely many common roots.
However, the polynomial compatibility equations
generally have  large number of terms,
and this method is currently out of our reach.

We take the opposite direction.
By successively differentiating the web curvature,
we show that the higher order derivatives of web  curvature
satisfy a simple functional relation.

Let $\, K$ be the web curvature, \eqref{webcurvature}.
Define $\, K_{-1}, \, K_{0}$ by
\beq
d K \equiv K_{-1} \, \omega+ K_0 \, \theta, \mod \; \, K.    \n
\eeq
Set the ideal of functions
\begin{align}
\mathcal{K}^{-1} &= 0, \n \\
\mathcal{K}^{0} &= \langle \, K \, \rangle, \n \\
\mathcal{K}^{1} &= \langle \, K, \, K_{-1}, \, K_0  \, \rangle, \n \\
                             &= \mathcal{K}^0  +  \nabla \mathcal{K}^0, \n
\end{align}
where $\, \langle \, K \, \rangle$ represents the ideal of functions
generated by $\, K$, and
$\,  \mathcal{K}^0  +  \nabla \mathcal{K}^0$ represents the ideal of functions
generated by
$\, \mathcal{K}^0$, and  the derivatives of  the elements in  $\, \mathcal{K}^0$.
Inductively define the curvature ideals
$\, \mathcal{K}^s$ for $\, s=0, \, 1, \, 2, \, . . . \, $  by
\beq
\mathcal{K}^{s+1} = \mathcal{K}^s  +   \nabla \mathcal{K}^s. \n
\eeq

Consider the local ring $\, \mathcal{O}(\textnormal{x$_0$})$
on  an infinitesimal neighborhood
of a reference point x$_0 \in M$.
If  $\, \mathcal{K}^{s_0}$ contains a unit in the local ring
$\, \mathcal{O}(\textnormal{x$_0$})$, i.e., a function nonzero at x$_0$,
then $\, \mathcal{K}^s = \mathcal{K}^{s_0} = \mathcal{O}(\textnormal{x$_0$})$
for $\, s \geq s_0$.
The sequence of  ideals $\,\{ \,  \K^s \, \}$   have  a meaning
only when the web curvature vanishes to certain order at  x$_0$.

Gronwall conjecture is equivalent to the equation
\beq
\K^0 /  \K^{-1} = 0. \n
\eeq
The main result of  this section
is a weaker  version of this equation;
\beq\label{Kmain}
\K^4 / \K^3 = 0.
\eeq

Recall that a  pencil   is a foliation on $\, \PP^2$
defined by a line in $\, (\PP^2)^*$.
\begin{theorem}\label{main}
Let $\, \W$ be a linear 3-web on a connected open subset $\, M \subset \PP^2$
which admits another  distinct linearization.
Let $\, \textnormal{x}_0 \in M$ be a reference point.
Suppose the web curvature of   $\, \W$ vanishes at  $\, \textnormal{x}_0$  up to order three.
Then the web curvature vanishes identically,
and $\, \W$ is algebraic.

In case $\, \W$ contains a pencil,
the same result holds with
the web curvature vanishing at  $\, \textnormal{x}_0$ up to order two.

In case $\, \W$ contains two pencils,
the same result holds with
the web curvature vanishing at  $\, \textnormal{x}_0$  up to order one.
\end{theorem}
\noi
The general case follows from \eqref{Kmain}
by the uniqueness theorem of ODE, see Appendix.
The cases containing pencils follow  similarly
from \eqref{Kmain1}, and  \eqref{Kmain2}.

Theorem \ref{main} is far from the   proof of Gronwall conjecture.
It is not likely  that
one could deduce from
the condition on the curvature ideal \eqref{Kmain} alone
that the web curvature must vanish.
The full proof would require a further differential analysis.

It is  nevertheless evident  that
a generic linear 3-web does admit  a unique linearization.
Theorem \ref{main} provides an explicit criterion
for this unique linearization
among  the class of linear 3-webs
with the web curvature that vanishes to certain order at a point.

\begin{example}
Let $\, (x, y)$ be the standard coordinate  of $\, \C^2$.
Fix a reference point $ \,(0, \, 0) \in \C^2$.
Consider a linear 3-web $\, \W$ defined by
the following three 1-forms
in a neighborhood of  $\, (0, \, 0)$.
\begin{align}\label{example3form}
\omega^1&=    d y,  \\
\omega^2&=  - h(x, \, y) \,  d x, \n \\
\omega^3&=  - ( dy - h(x, \, y) \,  d x), \n
\end{align}
where $\, h(0, 0) \ne 0$.
It is easily checked that $\, \W$ is linear when
$\, h(x, y)$ satisfies Burgers' equation
\beq\label{Burgers}
h_x + h \, h_y = 0.
\eeq
A short computation shows that the web curvature is given by
\begin{align}
K & = \frac{h_x \, h_y - h \, h_{xy}}{h^3}, \n \\
    &= \frac{h_{yy}}{h}. \n
\end{align}

Let  $\, h(x,y)$ be  the solution of  the following
initial value problem for Burgers' equation.
\begin{align}
h(0, y) & = 1 + y^{s+3} \, g(y), \n \\
h_x + h \, h_y &= 0, \n
\end{align}
for an integer $\, s \geq 1$, and an arbitrary nonzero analytic function $\, g(y)$ on $\, y$-axis.
By Cauchy-Kovalevsky theorem,
the analytic initial data can be uniquely thickened
to a solution in a neighborhood of $\, y$-axis.

Consider  a small neighborhood  $\, M$ of the origin $(0,0)$
on which $\, h(x,y)$ is nonzero, so that the three 1-forms \eqref{example3form}
define a linear 3-web.
\eqref{Burgers} implies that
every $\, x$-derivative of $\, h(x,y)$ can be replaced by
$\, y$-derivative up to scaling by nonzero terms, and modulo lower order terms.
From the formula for the web curvature,
it can be verified that the sequence of curvature ideals are  generated by
\begin{align}
\K^0 &= \langle \, \frac{\partial^2 h}{\partial y^2} \,   \rangle,  \n \\
\K^1 &= \langle \,\frac{\partial^2 h}{\partial y^2}, \, \frac{\partial^3 h}{\partial y^3} \,   \rangle,  \n \\
&... \n \\
\K^s &= \langle \, \frac{\partial^2 h}{\partial y^2}, \, \frac{\partial^3 h}{\partial y^3}, \, ... \,
                              \frac{\partial^{s+2} h}{\partial y^{s+2}} \,  \, \rangle. \n
\end{align}

The initial data satisfied by $\, h(x,y)$ implies that
the web curvature of the 3-web $\, \W$
vanishes   to order at least $\, s$ at   $\, (0,0)$.
The web curvature certainly does not vanish identically.
Hence $\, \W$ admits a unique local linearization for $\, s \geq 1$.
\end{example}

More generally,
a linear 3-web on an open subset $\, M \subset \PP^2$
is defined by a set of three disjoint curves in $\, M^* \subset (\PP^2)^*$.
Since a curve on a surface is locally described by a single function,
there exist roughly three arbitrary functions of one variable worth
local linear 3-webs.
Imposing the condition that the  web curvature vanishes
to a finite order  at a single point,
we get a subset of  finite codimension
in the moduli space of linear 3-webs.
Theorem \ref{main} thus provides a criterion for the unique linearization
for a fairly large subset of linear 3-webs.

In the following three subsections,
we present the proof of Theorem \ref{main}
for each case.
For computational simplicity,
scale $\, T_1 =1$,
and set accordingly
$\, \phi_{1,1} = \frac{-4 \,T_{{1,-1}}+ 1-  \,A_{{5}}}{9} \omega +A_{{4}}\theta$.

%-------------------------------------------------
\subsection{General case}\label{sec41}
Theorem \ref{main} is equivalent to
that the curvature ideal $\, \K^3$ is differentially closed;
\beq
d  \K^3 \equiv 0, \mod \; \, \K^3. \n
\eeq

Note  that one can solve for $\, A_{5,0}$ from  \eqref{webcurvature}.
Differentiating the web curvature \eqref{webcurvature}, one gets
\beq
\K^1 = \K^0  +  \, \langle \, K_1, \, K_2 \, \rangle, \n
\eeq
where
\begin{align}\label{K12}
K_1&=-T_{{1,-1}}A_{{4}}+7\,A_{{5}}A_{{4}}+5\,A_{{5}}T_{{1,-1}}A_{{4}}+15\,A
_{{4,0}}A_{{4}}-B_{{2}}+5\,A_{{4}}{T_{{1,-1}}}^{2}+2\,A_{{4}}
-108\,{A_{{4}}}^{3}  \\
&\quad +\frac{9}{2}\,T_{{0,0}}A_{{5}}  -3\,A_{{4}}T_{{1,-1,-1}}, \n \\
K_2&=-72\,{A_{{4}}}^{2}A_{{5}}+8\,A_{{4,0}}A_{{5}}-3\,B_{{1}}-B_{{3}}-5\,A_
{{4,0}}-4\,T_{{1,-1}}A_{{4,0}}
\n.
\end{align}
Note that one can solve for $\, \{ \, B_1, \, B_2 \, \}$ from $\, \{ K_1, \, K_2 \, \}$.

Differentiating \eqref{K12},  one gets
\beq
\K^2 = \K^1 +  \, \langle \, K_3, \, K_4 \, \rangle, \n
\eeq
where
\begin{align}\label{K34}
K_3&=-189\,{A_{{4}}}^{2}A_{{4,0}}-54\,T_{{0,0}}A_{{4}}-\frac{1}{2}\,A_{{4,0}}A_{{5}
}-T_{{1,-1}}A_{{4,0}}-36\,{A_{{4}}}^{2}T_{{1,-1}}+36\,{A_{{4}}}^{2}A_{
{5}}-B_{{5}}   \\
&\; \; -108\,{A_{{5}}}^{2}{A_{{4}}}^{2}+45\,{T_{{1,-1}}}^{2}{A_{{
4}}}^{2}-{\frac {81}{4}}\,{T_{{0,0}}}^{2}+648\,{A_{{4}}}^{4}-27\,T_{{0
,0}}T_{{1,-1}}A_{{4}}+{\frac {27}{2}}\,T_{{0,0}}A_{{5}}A_{{4}} \n \\
&\; \; -27\,{A_{{4}}}^{2}T_{{1,-1,-1}}+12\,A_{{4,0}}{A_{{5}}}^{2}-\frac{3}{2}\,A_{{5}}B_{{3}}
+45\,T_{{1,-1}}{A_{{4}}}^{2}A_{{5}}-A_{{4,0}}T_{{1,-1}}A_{{5}}+2\,A_{{
4,0}}-36\,{A_{{4}}}^{2} \n \\
&\; \; +15\,{A_{{4,0}}}^{2}+5\,A_{{4,0}}{T_{{1,-1}}}^{2}-3\,A_{{4,0}}T_{{1,-1,-1}},\n \\
K_4&=-3\,T_{{0,0}}A_{{4,0}}+5\,A_{{4,0}}A_{{4}}-8\,A_{{4}}A_{{4,0}}A_{{5}}+
4\,A_{{4}}A_{{4,0}}T_{{1,-1}}+\frac{1}{2}\,B_{{4}}+72\,{A_{{4}}}^{3}A_{{5}}+A_
{{4}}B_{{3}} \n \\
&\;\; +\frac{3}{2}\,T_{{1,-1,-1}}T_{{0,0}}-\frac{5}{2}\,T_{{1,-1}}A_{{5}}T_{{0,0
}}-\frac{5}{2}\,{T_{{1,-1}}}^{2}T_{{0,0}}+54\,{A_{{4}}}^{2}T_{{0,0}}+T_{{0,0}}
A_{{5}}+\frac{1}{2}\,T_{{0,0}}T_{{1,-1}} \n \\
&\;\;-T_{{0,0}}. \n
\end{align}
Note that one can solve for $\, \{ \, B_4, \, B_5 \, \}$ from $\, \{ K_3, \, K_4 \, \}$.

Differentiating \eqref{K34},  one gets
\beq
\K^3 = \K^2 +  \, \langle \, K_5 \, \rangle, \n
\eeq
where
\begin{align}\label{K5}
K_5&=
-72\,A_{{4,0}}T_{{1,-1}}{A_{{4}}}^{2}+72\,A_{{4,0}}T_{{0,0}}A_{{4}}-{
\frac {27}{2}}\,{T_{{0,0}}}^{2}A_{{5}}-90\,{A_{{4}}}^{2}A_{{4,0}}-3\,T
_{{0,0}}A_{{4}}-\frac{13}{3} \,A_{{4,0}}A_{{5}} \\
&\; \; +\frac{1}{2}\,T_{{1,-1}}A_{{4,0}}+24\,{A_{{4}}}^{2}A_{{5}}+\frac{1}{2}\,B_{{6}}
-24\,{A_{{5}}}^{2}{A_{{4}}}^{2}-1296\,{A_{{4}}}^{4}A_{{5}}+\frac{1}{3}\,B_{{3}}
-{\frac {27}{4}}\,{T_{{0,0}}}^{2}\n \\
&\;\; +\frac{3}{2}\,T_{{0,0}}T_{{1,-1}}A_{{4}}
-{\frac {75}{2}}\,T_{{0,0}}A_{{5}}A_{{4}}+
2\,A_{{4,0}}B_{{3}}-18\,{A_{{4}}}^{2}B_{{3}}+8\,{A_{{4,0}}}^{2}T_{{1,-
1}}+60\,T_{{1,-1}}{A_{{5}}}^{2}{A_{{4}}}^{2}\n \\
&\;\;+60\,{T_{{1,-1}}}^{2}A_{{5}}{A_{{4}}}^{2}
 +4\,T_{{1,-1,-1}}A_{{4,0}}A_{{5}}-2\,T_{{1,-1,-1}}A_{{4,0}}T_{{1,-1}}-16\,{A_{{4,0}}}^{2}A_{{5}}
 +\frac{8}{3}\,A_{{4,0}}{A_{{5}}}^{2}-\frac{1}{3}\,A_{{5}}B_{{3}} \n \\
 &\;\; +\frac{10}{3}\,A_{{4,0}}{T_{{1,-1}}}^{3}-36\,T_{{1,-1,-1}}{
A_{{4}}}^{2}A_{{5}}-\frac{10}{3}\,A_{{4,0}}{T_{{1,-1}}}^{2}A_{{5}}
-\frac{15}{2}\,T_{{0,0}}{T_{{1,-1}}}^{2}A_{{4}} \n \\
&\;\; -\frac{15}{2}\,T_{{0,0}}T_{{1,-1}}A_{{5}}A_{{4}}+{
\frac {27}{4}}\,{T_{{0,0}}}^{2}T_{{1,-1}}-\frac{1}{6}\,B_{{3}}T_{{1,-1}}-\frac{1}{2}\,
T_{{1,-1,-1}}B_{{3}}+\frac{5}{6}\,{T_{{1,-1}}}^{2}B_{{3}} \n \\
&\;\;+288\,A_{{4,0}}A_{{5}
}{A_{{4}}}^{2}-648\,T_{{0,0}}{A_{{4}}}^{3}-12\,T_{{1,-1}}{A_{{4}}}^{2}
A_{{5}}+\frac{9}{2}\,A_{{4}}T_{{1,-1,-1}}T_{{0,0}}-{\frac {20}{3}}\,A_{{4,0}}{A_{{5}}}^{2}T_{{1,-1}} \n \\
&\;\;+\frac{5}{6}\,A_{{5}}B_{{3}}T_{{1,-1}}+{\frac {25}{6}}\,
A_{{4,0}}T_{{1,-1}}A_{{5}}+\frac{5}{3}\,A_{{4,0}}+10\,{A_{{4,0}}}^{2}+\frac{7}{2}\,A_{
{4,0}}{T_{{1,-1}}}^{2}-\frac{5}{2}\,A_{{4,0}}T_{{1,-1,-1}}.
 \n
\end{align}
Note that one can solve for $\, \{ \, B_6 \, \}$ from $\, \{ \, K_5 \, \}$.

Differentiating \eqref{K5},  one finally  gets
\beq
d  K_5 \equiv 0, \mod \; \, \K^3. \n
\eeq

\subsection{Case with one pencil}\label{sec42}
Assume  $\, \sigma_3$-foliation is a pencil,  \eqref{sigmaformula}.
From the structure equation,
\begin{align}
d Z^2 &\equiv - Z^0 \, \theta,  \mod \; \, Z^2, \n \\
d Z^0  &\equiv - T_{0,0} \, Z^1 \, \theta,  \mod \; \, Z^2, \, Z^0. \n
\end{align}
$\, \sigma_3$-foliation is a pencil when $\, T_{0,0} = 0$.
Successive derivatives of the equation  $\, T_{0,0} = 0$ then imply that
$\, B_1, \, B_4$, and  $\, B_6$ are zero.

Theorem \ref{main} is equivalent to
that the curvature ideal $\, \K^2$ is differentially closed;
\beq\label{Kmain1}
d  \K^2 \equiv 0, \mod \; \, \K^2.
\eeq
Differentiating the web curvature \eqref{webcurvature}, one gets
\beq
\K^1 = \K^0  +  \, \langle \, K_1, \, K_2 \, \rangle, \n
\eeq
where
\begin{align}\label{KK12}
K_1&=-T_{{1,-1}}A_{{4}}+5\,A_{{4}}{T_{{1,-1}}}^{2}+7\,A_{{5}}A_{{4}}+15\,A_
{{4,0}}A_{{4}}-108\,{A_{{4}}}^{3}+5\,A_{{5}}T_{{1,-1}}A_{{4}}   \\
&\quad -3\,A_{{4 }}T_{{1,-1,-1}}+2\,A_{{4}}-B_{{2}},  \n \\
K_2&=-4\,T_{{1,-1}}A_{{4,0}}-5\,A_{{4,0}}+8\,A_{{4,0}}A_{{5}}-B_{{3}}
-72\,{A_{{4}}}^{2}A_{{5}}.\n
\end{align}
Note that one can solve for $\, \{ \, B_2, \, B_3 \, \}$ from $\, \{ K_1, \, K_2 \, \}$.

Differentiating \eqref{KK12},  one gets
\beq
\K^2 = \K^1 +  \, \langle \, K_3  \, \rangle, \n
\eeq
where
\begin{align}\label{KK3}
K_3&=45\,{T_{{1,-1}}}^{2}{A_{{4}}}^{2}+45\,T_{{1,-1}}{A_{{4}}}^{2}A_{{5}}+5
\,A_{{4,0}}T_{{1,-1}}A_{{5}}-36\,{A_{{4}}}^{2}T_{{1,-1}}-36\,{A_{{4}}}
^{2}+36\,{A_{{4}}}^{2}A_{{5}}   \\
&\quad -T_{{1,-1}}A_{{4,0}}+648\,{A_{{4}}}^{4}+5
\,A_{{4,0}}{T_{{1,-1}}}^{2}-27\,{A_{{4}}}^{2}T_{{1,-1,-1}}+15\,{A_{{4,0
}}}^{2}-3\,A_{{4,0}}T_{{1,-1,-1}} \n \\
&\quad -189\,{A_{{4}}}^{2}A_{{4,0}}+7\,A_{{4,0}}A_{{5}}-B_{{5}}+2\,A_{{4,0}}.
\n
\end{align}
Note that one can solve for $\, \{ \,   B_5 \, \}$ from $\, \{ \, K_3  \, \}$.

Differentiating \eqref{KK3},  one finally  gets
\beq
d  K_3 \equiv 0, \mod \; \, \K^2. \n
\eeq

\subsection{Case with two pencils}\label{sec43}
Assume $\, \sigma_1$, and $\, \sigma_2$-foliations are pencils, \eqref{sigmaformula}, i.e.,
\begin{align}
d \sigma_i   &\equiv \sigma_i' \, \alpha_i, \mod \; \, \sigma_i,  \n \\
d \sigma_i'  &\equiv 0,  \quad \;\,    \mod \; \, \sigma_i, \, \sigma_i', \n
\end{align}
for $\, i=1, \, 2\,$.
From the structure equation,
a computation shows that this  implies
\begin{align}\label{XXX12}
 A_{{5,0}}&=-5\,A_{{4}}-T_{{0,0}}+8\,A_{{5}}A_{{4}}-4\,T_{{1,-1}}A_{{4}},   \\
T_{{1,-1,-1}}&=-\frac{1}{3}\,{A_{{5}}}^{2}+T_{{1,-1}}A_{{5}}-A_{{4,0}}
+\frac{1}{3}+9\,{A_{{4}}}^{2}+\frac{4}{3}\,{T_{{1,-1}}}^{2}+\frac{1}{3}\,T_{{1,-1}}. \n
\end{align}
In this case,  the web curvature is given by  $\, K = - T_{0,0}$.

Theorem \ref{main} is   equivalent to
that the curvature ideal $\, \K^1$ is differentially closed;
\beq\label{Kmain2}
d \K^1 \equiv 0, \mod \; \, \K^1.
\eeq
From \eqref{T00},   we have
\beq
\K^1 = \K^0   +   \, \langle \, K_1 \, \rangle, \n
\eeq
where
\beq\label{KKK1}
K_1 = B_1.\n
\eeq
It suffices to  show that
\beq
d B_1  \equiv 0, \mod \; \, T_{0,0}, \,B_1. \n
\eeq

Differentiating   \eqref{XXX12}, one gets a set of two equations,
which give
\begin{align}\label{XXX34}
B_{{2}}&=-2\,T_{{1,-1}}A_{{4}}-135\,{A_{{4}}}^{3}+2\,A_{{5}}T_{{1,-1}}A_{{4}}
+\frac{7}{6}\,T_{{0,0}}A_{{5}}+18\,A_{{4,0}}A_{{4}}+\frac{2}{3}\,T_{{0,0}}T_{{1,-1}}   \\
&\quad +7\,A_{{5}}A_{{4}}+\frac{4}{3}\,T_{{0,0}}+A_{{4}}{T_{{1,-1}}}^{2}
+{A_{{5}}}^{2}A_{{4}}+A_{{4}},  \n \\
B_{{3}}&=-B_{{1}}-4\,T_{{1,-1}}A_{{4,0}}-5\,A_{{4,0}}-6\,T_{{0,0}}A_{{4}}
 -72\,{A_{{4}}}^{2}A_{{5}}+8\,A_{{4,0}}A_{{5}}. \n
\end{align}

Differentiating   \eqref{XXX34}, one gets another  set of two equations,
which give
\begin{align}\label{XXX56}
B_{{4}}&=18\,B_{{1}}A_{{4}}+{\frac {14}{9}}\,T_{{1,-1}}A_{{5}}T_{{0,0}}
+\frac{1}{9}\,{T_{{1,-1}}}^{2}T_{{0,0}}+25\,T_{{0,0}}A_{{4,0}}-{\frac {35}{9
}}\,T_{{0,0}}+{\frac {13}{9}}\,T_{{0,0}}{A_{{5}}}^{2}  -387\,{A_{{4}}}^{2}T_{{0,0}}    \\
&\quad-{\frac {74}{9}}\,T_{{0,0}}T_{{1,-1}}
+{\frac {58}{9}}\,T_{{0,0}}A_{{5}}, \n \\
B_{{5}}&=\frac{4}{3}\,B_{{1}}+2\,A_{{4,0}}T_{{1,-1}}A_{{5}}+18\,{A_{{4
,0}}}^{2}-189\,{A_{{4}}}^{2}A_{{4,0}}+A_{{4,0}}{T_{{1,-1}}}^{2}+A_{{4,0
}}{A_{{5}}}^{2}+\frac{7}{6}\,A_{{5}}B_{{1}} \n  \\
&\quad +\frac{2}{3}\,B_{{1}}T_{{1,-1}}+18\,T_{{1,-
1}}{A_{{4}}}^{2}A_{{5}}-\frac{7}{2}\,T_{{0,0}}T_{{1,-1}}A_{{4}}+T_{{0,0}}A_{{5
}}A_{{4}}+405\,{A_{{4}}}^{4}+9\,{A_{{5}}}^{2}{A_{{4}}}^{2}-{\frac {13}
{4}}\,{T_{{0,0}}}^{2} \n \\
&\quad -{\frac {35}{2}}\,T_{{0,0}}A_{{4}}+9\,{T_{{1,-1}}
}^{2}{A_{{4}}}^{2}-45\,{A_{{4}}}^{2}+A_{{4,0}}+7\,A_{{4,0}}A_{{5}}+36
\,{A_{{4}}}^{2}A_{{5}}-45\,{A_{{4}}}^{2}T_{{1,-1}}-2\,T_{{1,-1}}A_{{4,0}}.  \n
\end{align}
Note $\,  B_4 \equiv 0,  \mod \; \,  T_{0,0}, \, B_1$.
From  \eqref{B1},
this implies $\, d B_1 \equiv 0, \mod \; \,  T_{0,0}, \, B_1$.

Differentiating   \eqref{XXX56}, one again gets a set of two equations,
which give
\begin{align}
 B_{{6}}&=43\,B_{{1}}A_{{4,0}}-1053\,T_{{0,0}}{A_{{4}}}^{3}-171
\,B_{{1}}{A_{{4}}}^{2}-{\frac {382}{3}}\,T_{{0,0}}T_{{1,-1}}A_{{4}}-{
\frac {74}{9}}\,B_{{1}}T_{{1,-1}}-{\frac {35}{9}}\,B_{{1}} \n  \\
&\quad +{\frac {65}{3}}\,T_{{0,0}}{T_{{1,-1}}}^{2}A_{{4}}+{\frac {127}{6}}\,{T_{{0,0}}}^{
2}A_{{5}}+{\frac {38}{3}}\,{T_{{0,0}}}^{2}T_{{1,-1}}+{\frac {166}{3}}
\,T_{{0,0}}T_{{1,-1}}A_{{5}}A_{{4}}+{\frac {64}{3}}\,{T_{{0,0}}}^{2} \n \\
&\quad +{\frac {101}{3}}\,T_{{0,0}}{A_{{5}}}^{2}A_{{4}}+{\frac {545}{3}}\,T_{{0
,0}}A_{{5}}A_{{4}}-{\frac {115}{3}}\,T_{{0,0}}A_{{4}}+{\frac {58}{9}}
\,A_{{5}}B_{{1}}+\frac{1}{9}\,B_{{1}}{T_{{1,-1}}}^{2}+{\frac {13}{9}}\,B_{{1}}
{A_{{5}}}^{2}\n \\
&\quad +{\frac {14}{9}}\,B_{{1}}T_{{1,-1}}A_{{5}}-99\,A_{{4,0}}T_{{0,0}}A_{{4}}, \n
\end{align}
and
\beq\label{XXX7}
A_{{4}} \left( T_{{1,-1}}+1 \right) B_{{1}} \equiv 0, \mod \; \,T_{0,0}.
\eeq
Note that from \eqref{XXX7},
if  $\, A_4 (T_{1,-1}+1) \ne 0$ at the reference point x$_0$,
we get an additional equation $\, B_1 \equiv 0, \mod \;  T_{0,0}$.
In this case,  the vanishing of the web curvature to order zero at x$_0$ suffices.

%-------------------------------------------------
%-------------------------------------------------
\section{Concluding remark}\label{sec5}
The differential relation satisfied by the derivatives of web curvature
suggests an approach
toward the full proof of Gronwall conjecture.

Let $\, \vec{K} = (K, \, K_1, \, K_2, \, ... \, K_l)^t$
denote the derivatives of web curvature $\, K$
that satisfy    \eqref{Kmain}, \eqref{Kmain1}, or \eqref{Kmain2}
($\, l =5, \, 3, \, 1$  respectively).
The differential relation implies that
\beq\label{vecK}
d \vec{K} = \gamma \, \vec{K},
\eeq
for an $\,(l+1)$-by-$\, (l+1)$ matrix  1-form $\, \gamma$.
Set $\, d \gamma - \gamma \w \gamma =  Q \, \omega \w \theta$.
Differentiating \eqref{vecK}, one gets
\beq\label{QK}
Q \, \vec{K} = 0.
\eeq
Gronwall conjecture would follow by showing that
$\, det(Q)=0$ and  \eqref{QK}
imply  $\, K =0$.

Consider the simplest two pencil case with $\, A_4 (T_{1,-1}+1) \ne 0$.
There exists a 1-form $\, \gamma$ such that
\beq
d K = \gamma \, K. \n
\eeq
Let $\, d \gamma = Q \, \omega \w \theta$,
and
let $\, \vec{Q} = (Q, \, Q_1, \, Q_2, \, ... ) \, $ denote  the successive derivatives of $\, Q$.
If one can show that there are at most finitely many
common roots for $\, \vec{Q}$,
a short analysis  implies   that $\, K = 0$.

$\, Q$ is a degree 6  polynomial in five variables.
The degree  and the size of $\, Q_i$  increase as one differentiates.
The analysis of the root structure of    $\, \vec{Q}$,
although it may not be feasible manually,
would be the first step to prove, or disprove, Gronwall conjecture.
Bol gave a case by case analysis of a set of  quasi-algebraic linear 3-webs
in  an attempt to  find a counterexample to the conjecture, \cite{Bol}.

\renewcommand{\theequation}{A-\arabic{equation}}
 %redefine the command that creates the equation no.
\setcounter{equation}{0}  % reset counter
\section*{Appendix}

\noi
A-1. \textbf{Uniqueness theorem of ODE}

Let $\, U \subset \C^n$ be a connected open subset containing the origin $\, 0 \in \C^n$.
Let $\, \vec{f} = (f^1, \, f^2, \, ... \, f^m)^t$ denote a $\, \C^m$-valued function on $\, U$.
Let $\, \gamma$ be an $m$-by-$m$ matrix   1-form on $\, U$.

Consider the   initial value problem for the linear ODE;
\begin{align}
d \vec{f} &= \gamma \, \vec{f}, \n \\
   \vec{f}(0)&= \vec{f}_0. \n
\end{align}
It admits a unique solution,  if the solution exists.
In particular, if $\,  \vec{f}_0=0$,
$\,  \vec{f} \equiv 0$ is the only solution.

\vsp{1pc}
\noi
A-2. \textbf{Formulae}

We collect the exact formulae of
some of the  long expressions omitted in the main text.

\vsp{.5pc}

\noi
$A_{4,0,-1}$:
\begin{align}
&-\frac{1}{12 \, T_1^2}  (-9\,T_{{0,0}}T_{{1}}A_{{4}}-9\,A_{{5,0}}T_{{1}}A_{{4}}+
72\,A_{{5}}{A_{{4}}}^{2}+72\,T_{{1,-1}}{A_{{4}}}^{2}+27\,{A_{{4}}}^{2}
{T_{{1}}}^{2}+32\,A_{{4,0}}{T_{{1}}}^{3}+6\,B_{{1}}{T_{{1}}}^{2}\n \\
&-20\,A_{{4,0}}T_{{1}}T_{{1,-1}}-20\,A_{{4,0}}T_{{1}}A_{{5}})  \n
\end{align}

%-------------------------------------------------------------------------------------
\noi
$A_{5,0,-1}$:
\begin{align}
&-\frac{1}{6 \, T_1^2}(   38\,A_{{5}}T_{{1,-1}}A_{{4}}+21\,T_{{1}}A_{{5}}T_{{0,0}}
+42\,A_{{4,0}}A_{{4}}T_{{1}}-6\,B_{{2}}{T_{{1}}}^{2}-18\,A_{{5,0}}T_{{
1}}A_{{5}}-12\,A_{{5,0}}T_{{1}}T_{{1,-1}}-216\,{A_{{4}}}^{3}\n \\
&-6\,T_{{0,0}}T_{{1}}T_{{1,-1}}-10\,A_{{4}}{T_{{1,-1}}}^{2}+48\,{A_{{5}}}^{2}A_{{4
}}+40\,{T_{{1}}}^{2}A_{{4}}A_{{5}}+6\,A_{{4}}T_{{1,-1,-1}}T_{{1}}+12\,
A_{{5,0}}{T_{{1}}}^{3}+2\,A_{{4}}{T_{{1}}}^{4}\n \\
&-3\,T_{{0,0}}{T_{{1}}}^{3}-4\,{T_{{1}}}^{2}T_{{1,-1}}A_{{4}} ) \n
\end{align}

%-------------------------------------------------------------------------------------
\noi
$T_{1,-1,-1,-1}$:
\begin{align}
&\frac{1}{36\, T_1^2} \,(-1161\,{A_{{4}}}^{2}{T_{{1}}}^{2}+147\,A_{{4,0}}{T_{{1}}
}^{3}-351\,T_{{0,0}}T_{{1}}A_{{4}}-189\,A_{{5,0}}T_{{1}}A_{{4}}+3456\,
A_{{5}}{A_{{4}}}^{2}-1620\,T_{{1,-1}}{A_{{4}}}^{2} \n \\
&+204\,A_{{4,0}}T_{{1}}T_{{1,-1}}-216\,A_{{4,0}}T_{{1}}A_{{5}}
+27\,B_{{3}}{T_{{1}}}^{2}+8\,{T_{{1}}}^{2}A_{{5}}T_{{1,-1}}
-16\,{T_{{1}}}^{4}A_{{5}}-160\,{T_{{1,-1}}}^{2}A_{{5}}\n \\
&+120\,T_{{1,-1,-1}}T_{{1}}A_{{5}}+45\,B_{{1}}{T_{{1}}}^{
2}+16\,{T_{{1}}}^{6}+216\,T_{{1,-1,-1}}T_{{1}}T_{{1,-1}}+24\,{T_{{1}}}
^{4}T_{{1,-1}}+48\,{T_{{1,-1}}}^{2}{T_{{1}}}^{2}\n \\
&-160\,{T_{{1,-1}}}^{3}-72\,T_{{1,-1,-1}}{T_{{1}}}^{3}) \n
\end{align}

%-------------------------------------------------------------------------------------
\noi
$T_{1,-1,-1,0}$:
\begin{align}
&\frac{1}{24\, T_1^2}\,(80\,{A_{{5}}}^{2}A_{{4}}+2\,{T_{{1}}}^{2}T_{{1,-1}}A_{{4}}
-648\,{A_{{4}}}^{3}+118\,A_{{5}}T_{{1,-1}}A_{{4}}-3\,T_{{1}}A_{{5}}T_{{0,0}}
+90\,A_{{4,0}}A_{{4}}T_{{1}}-10\,A_{{5,0}}T_{{1}}A_{{5}} \n \\
&
-16\,A_{{5,0}}T_{{1}}T_{{1,-1}}-120\,T_{{0,0}}T_{{1}}T_{{1,-1}}
-40\,{T_{{1}}}^{2}A_{{4}}A_{{5}}-18\,A_{{4}}T_{{1,-1,-1}}T_{{1}}
-6\,B_{{2}}{T_{{1}}}^{2}-34\,A_{{4}}{T_{{1,-1}}}^{2} \n \\
&
+4\,A_{{5,0}}{T_{{1}}}^{3}-40\,A_{{4}}{T_{{1}}}^{4}
-24\,T_{{0,0}}{T_{{1}}}^{3}) \n
\end{align}

%-------------------------------------------------------------------------------------
\noi
$B_{1,-1}$:
\begin{align}
&\frac{1}{4\, T_1^2} (12\,T_{{1,-1}}B_{{1}}T_{{1}}-8\,B_{{1}}{T_{{1}}}^{3}+8\,B
_{{1}}T_{{1}}A_{{5}}-24\,T_{{0,0}}T_{{1,-1}}A_{{4}}-3\,{T_{{0,0}}}^{2}
T_{{1}}+3\,T_{{0,0}}A_{{5,0}}T_{{1}}-24\,T_{{0,0}}A_{{5}}A_{{4}}\n \\
&
-21\,{T_{{1}}}^{2}A_{{4}}T_{{0,0}})\n
\end{align}

%-------------------------------------------------------------------------------------
\noi
$B_{2,-1}$:
\begin{align}
&\frac{1}{12\, T_1^3}   (-6\,B_{{4}}{T_{{1}}}^{3}-40\,B_{{2}}{T_{{1}}}^{4}+24\,A_
{{5,0}}{T_{{1}}}^{2}A_{{4,0}}-264\,A_{{4,0}}T_{{1}}A_{{5}}A_{{4}}+57\,
{T_{{1}}}^{3}A_{{4,0}}A_{{4}}+9\,A_{{4}}B_{{1}}{T_{{1}}}^{2}\n \\
&
+72\,A_{{4,0}}{T_{{1}}}^{2}T_{{0,0}}+9\,A_{{4}}B_{{3}}{T_{{1}}}^{2}
-84\,A_{{4,0}}T_{{1}}T_{{1,-1}}A_{{4}}-108\,A_{{5,0}}T_{{1}}{A_{{4}}}^{2}
+1512\,{A_{{4}}}^{3}A_{{5}}+864\,{A_{{4}}}^{3}T_{{1,-1}}\n \\
&
+432\,{T_{{1}}}^{2}{A_{{4}}}^{3}+28\,B_{{2}}{T_{{1}}}^{2}T_{{1,-1}}
+28\,B_{{2}}{T_{{1}}}^{2}A_{{5}})  \n
\end{align}

%-------------------------------------------------------------------------------------
\noi
$B_{3,-1}$:
\begin{align}
&\frac{1}{24\, T_1^3}\,(-152\,A_{{4,0}}T_{{1}}T_{{1,-1}}A_{{5}}+12\,A_{{5,0}}T_{
{1}}T_{{1,-1}}A_{{4}}-816\,A_{{5,0}}T_{{1}}A_{{5}}A_{{4}}-36\,T_{{1,-1
}}A_{{4}}T_{{0,0}}T_{{1}}+2421\,T_{{0,0}}T_{{1}}A_{{5}}A_{{4}}\n \\
&
-162\,A_{{4}}B_{{2}}{T_{{1}}}^{2}+3510\,A_{{4,0}}T_{{1}}{A_{{4}}}^{2}
+40\,A_{{4,0}}T_{{1}}{T_{{1,-1}}}^{2}+60\,{A_{{5,0}}}^{2}{T_{{1}}}^{2}+12\,B_{{
1}}{T_{{1}}}^{4}-64\,B_{{3}}{T_{{1}}}^{4}\n \\
&
-168\,{A_{{4,0}}}^{2}{T_{{1}}}^{2}+4416\,{A_{{5}}}^{2}{A_{{4}}}^{2}
+108\,{T_{{0,0}}}^{2}{T_{{1}}}^{2}+24\,B_{{5}}{T_{{1}}}^{3}
-318\,{T_{{1,-1}}}^{2}{A_{{4}}}^{2}+348\,{T_{{1}}}^{4}{A_{{4}}}^{2}\n \\
&
-8\,{T_{{1}}}^{5}A_{{4,0}}-17496\,{A_{{4}}}^{4}-192\,A_{{4,0}}T_{{1}}{A_{{5}}}^{2}
-24\,A_{{4,0}}{T_{{1}}}^{2}T_{{1,-1,-1}}-18\,A_{{5,0}}{T_{{1}}}^{2}T_{{0,0}}
-84\,B_{{1}}{T_{{1}}}^{2}A_{{5}}\n \\
&+24\,B_{{1}}{T_{{1}}}^{2}T_{{1,-1}}
+88\,B_{{3}}{T_{{1}}}^{2}A_{{5}}+64\,B_{{3}}{T_{{1}}}^{2}T_{{1,-1}}
+1938\,T_{{1,-1}}{A_{{4}}}^{2}A_{{5}}+162\,{A_{{4}}}^{2}T_{{1,-1,-1}}T_{{1}} \n \\
&+3390\,{T_{{1}}}^{2}A_{{5}}{A_{{4}}}^{2}-30\,{T_{{1}}}^{2}T_{{1,-1}}{A_{{4}}}^{2}
 +16\,{T_{{1}}}^{3}T_{{1,-1}}A_{{4,0}}-160\,{T_{{1}}}^{3}A_{{4,0}}A_{{5}}+198\,{T_{{1}}}
^{3}A_{{4}}T_{{0,0}}-12\,{T_{{1}}}^{3}A_{{4}}A_{{5,0}})\n
\end{align}

%-------------------------------------------------------------------------------------
\noi
$B_{3,0}$:
\begin{align}
&\frac{1}{18\, T_1^3} \,(-5712\,A_{{4,0}}T_{{1}}A_{{5}}A_{{4}}+1722\,A_{{4,0}}T_{
{1}}T_{{1,-1}}A_{{4}}+120\,T_{{1,-1}}A_{{5}}T_{{0,0}}T_{{1}}+4\,T_{{1,
-1}}A_{{5}}A_{{5,0}}T_{{1}}\n \\
&+480\,T_{{1,-1,-1}}T_{{1}}A_{{5}}A_{{4}}
-402\,T_{{1,-1,-1}}T_{{1}}T_{{1,-1}}A_{{4}}+470\,{T_{{1}}}^{2}A_{{5}}T_
{{1,-1}}A_{{4}}+282\,A_{{5,0}}{T_{{1}}}^{2}A_{{4,0}}\n \\
&+108\,A_{{4}}B_{{3}}{T_{{1}}}^{2}-3024\,A_{{5,0}}T_{{1}}{A_{{4}}}^{2}
 -126\,B_{{2}}{T_{{1}}}^{2}T_{{1,-1}}+144\,B_{{2}}{T_{{1}}}^{2}A_{{5}}-27\,B_{{4}}{T_{{1}}
}^{3}-90\,B_{{2}}{T_{{1}}}^{4} \n \\
&+51408\,{A_{{4}}}^{3}A_{{5}}-14040\,{A_{
{4}}}^{3}T_{{1,-1}}-7344\,{T_{{1}}}^{2}{A_{{4}}}^{3}
+718\,{T_{{1,-1}}}^{3}A_{{4}}+16\,{T_{{1}}}^{5}A_{{5,0}} \n \\
&+80\,{T_{{1}}}^{6}A_{{4}}-6\,{T_
{{1}}}^{5}T_{{0,0}}+112\,A_{{5,0}}T_{{1}}{T_{{1,-1}}}^{2}-60\,A_{{5,0}
}{T_{{1}}}^{2}T_{{1,-1,-1}}+201\,{T_{{1,-1}}}^{2}T_{{0,0}}T_{{1}}\n \\
&
-32\,T_{{1,-1}}{A_{{5}}}^{2}A_{{4}}-610\,{T_{{1,-1}}}^{2}A_{{5}}A_{{4}}-99
\,T_{{1,-1,-1}}{T_{{1}}}^{2}T_{{0,0}}-324\,{A_{{4}}}^{2}T_{{0,0}}T_{{1
}}-8\,{T_{{1}}}^{3}A_{{5,0}}T_{{1,-1}}\n \\
&+3\,{T_{{1}}}^{3}T_{{0,0}}T_{{1,-1}}
+474\,{T_{{1}}}^{2}A_{{4}}{T_{{1,-1}}}^{2}-16\,{T_{{1}}}^{2}{A_{{5
}}}^{2}A_{{4}}-118\,{T_{{1}}}^{4}A_{{4}}A_{{5}}-300\,{T_{{1}}}^{3}A_{{
4}}T_{{1,-1,-1}} \n \\
&+132\,{T_{{1}}}^{4}T_{{1,-1}}A_{{4}}
+2\,{T_{{1}}}^{3}A_{{5,0}}A_{{5}}+60\,{T_{{1}}}^{3}A_{{5}}T_{{0,0}}
+162\,A_{{4}}B_{{1}}{T_{{1}}}^{2}+36\,A_{{4,0}}{T_{{1}}}^{2}T_{{0,0}}
+1140\,{T_{{1}}}^{3}A_{{4,0}}A_{{4}}) \n
\end{align}

%-------------------------------------------------------------------------------------
\noi
$B_{4,-1}$:
\begin{align}
&-\frac{1}{12\, T_1^3} \,(-44\,B_{{4}}{T_{{1}}}^{2}T_{{1,-1}}-32\,B_{{4}}{T_{{1}}
}^{2}A_{{5}}+32\,B_{{4}}{T_{{1}}}^{4}+72\,T_{{0,0}}A_{{4,0}}T_{{1}}T_{
{1,-1}}+72\,T_{{0,0}}A_{{4,0}}T_{{1}}A_{{5}}\n \\
&
-351\,A_{{4}}{T_{{0,0}}}^{2}T_{{1}}-1269\,{A_{{4}}}^{2}T_{{0,0}}{T_{{1}}}^{2}
-63\,B_{{1}}{T_{{1}}}^{2}T_{{0,0}}-24\,B_{{1}}{T_{{1}}}^{2}A_{{5,0}}
-9\,T_{{0,0}}B_{{3}}{T_{{1}}}^{2}\n \\
&
+63\,A_{{4,0}}{T_{{1}}}^{3}T_{{0,0}}-24\,B_{{1}}T_{{1}}T_{
{1,-1}}A_{{4}}+192\,B_{{1}}T_{{1}}A_{{5}}A_{{4}}-1080\,{A_{{4}}}^{2}T_
{{0,0}}T_{{1,-1}}+135\,A_{{5,0}}T_{{1}}A_{{4}}T_{{0,0}}\n \\
&
-1728\,{A_{{4}}}^{2}T_{{0,0}}A_{{5}}-48\,B_{{1}}{T_{{1}}}^{3}A_{{4}})\n
\end{align}

%-------------------------------------------------------------------------------------
\noi
$B_{5,-1}$:
\small{
\begin{align}
&\frac{1}{24\, T_1^4} \,(-43848\,{A_{{4}}}^{4}T_{{1,-1}}+718\,{T_{{1,-1}}}^{3}{A_
{{4}}}^{2}-32832\,{A_{{4}}}^{4}A_{{5}}-28080\,{T_{{1}}}^{2}{A_{{4}}}^{
4}+80\,{T_{{1}}}^{6}{A_{{4}}}^{2}\n \\
&
+114\,{T_{{1}}}^{4}{A_{{4,0}}}^{2}-96
\,{T_{{1}}}^{5}B_{{5}}+72\,B_{{5}}{T_{{1}}}^{3}T_{{1,-1}}+72\,B_{{5}}{
T_{{1}}}^{3}A_{{5}}-168\,{A_{{4,0}}}^{2}{T_{{1}}}^{2}T_{{1,-1}}\n \\
&
-32\,T_{{1,-1}}{A_{{5}}}^{2}{A_{{4}}}^{2}-610\,{T_{{1,-1}}}^{2}A_{{5}}{A_{{4}
}}^{2}+2160\,A_{{5,0}}T_{{1}}{A_{{4}}}^{3}+54\,{A_{{4}}}^{2}B_{{1}}{T_
{{1}}}^{2}\n \\
&-216\,{A_{{4}}}^{2}B_{{3}}{T_{{1}}}^{2}-6156\,{A_{{4}}}^{3}T
_{{0,0}}T_{{1}}+66\,A_{{4,0}}{T_{{1}}}^{3}B_{{3}}+306\,B_{{2}}{T_{{1}}
}^{3}T_{{0,0}}-9\,A_{{4}}B_{{4}}{T_{{1}}}^{3}\n \\
&+162\,A_{{4,0}}B_{{1}}{T_
{{1}}}^{3}+90\,B_{{2}}{T_{{1}}}^{3}A_{{5,0}}+2940\,{T_{{1}}}^{3}A_{{4,0
}}{A_{{4}}}^{2}-16\,{T_{{1}}}^{2}{A_{{5}}}^{2}{A_{{4}}}^{2}\n \\
&
-300\,{T_{{1}}}^{3}{A_{{4}}}^{2}T_{{1,-1,-1}}-6\,{T_{{1}}}^{5}A_{{4}}T_{{0,0}}+16
\,{T_{{1}}}^{5}A_{{4}}A_{{5,0}}-118\,{T_{{1}}}^{4}A_{{5}}{A_{{4}}}^{2}
-60\,A_{{4}}A_{{5,0}}{T_{{1}}}^{2}T_{{1,-1,-1}}\n \\
&
+201\,A_{{4}}{T_{{1,-1}}}^{2}T_{{0,0}}T_{{1}}-99\,A_{{4}}T_{{1,-1,-1}}{T_{{1}}}^{2}T_{{0,0}}-
342\,A_{{4}}A_{{4,0}}{T_{{1}}}^{2}T_{{0,0}}-402\,T_{{1,-1,-1}}T_{{1}}T
_{{1,-1}}{A_{{4}}}^{2}\n \\
&
+112\,A_{{4}}A_{{5,0}}T_{{1}}{T_{{1,-1}}}^{2}-
492\,A_{{4,0}}{T_{{1}}}^{2}A_{{4}}A_{{5,0}}+10416\,A_{{4,0}}T_{{1}}A_{
{5}}{A_{{4}}}^{2}+7050\,A_{{4,0}}T_{{1}}T_{{1,-1}}{A_{{4}}}^{2}\n \\
&
+480\,T_{{1,-1,-1}}T_{{1}}A_{{5}}{A_{{4}}}^{2}-720\,B_{{2}}{T_{{1}}}^{2}A_{{5
}}A_{{4}}-126\,B_{{2}}{T_{{1}}}^{2}T_{{1,-1}}A_{{4}}-8\,{T_{{1}}}^{3}A
_{{5,0}}T_{{1,-1}}A_{{4}}\n \\
&
+2\,{T_{{1}}}^{3}A_{{5,0}}A_{{5}}A_{{4}}+470
\,{T_{{1}}}^{2}T_{{1,-1}}{A_{{4}}}^{2}A_{{5}}+3\,{T_{{1}}}^{3}T_{{1,-1
}}A_{{4}}T_{{0,0}}+60\,{T_{{1}}}^{3}T_{{0,0}}A_{{5}}A_{{4}}\n \\
&
+4\,A_{{4}}T_{{1,-1}}A_{{5}}A_{{5,0}}T_{{1}}+120\,A_{{4}}T_{{1,-1}}A_{{5}}T_{{0,0
}}T_{{1}}+132\,{T_{{1}}}^{4}T_{{1,-1}}{A_{{4}}}^{2}-528\,{A_{{4,0}}}^{
2}{T_{{1}}}^{2}A_{{5}}\n \\
&
+474\,{T_{{1}}}^{2}{T_{{1,-1}}}^{2}{A_{{4}}}^{2}
-12\,B_{{6}}{T_{{1}}}^{4}+306\,{T_{{1}}}^{4}A_{{4}}B_{{2}}) \n
\end{align}
}

%-------------------------------------------------------------------------------------
\noi
$B_{5,0}$:
\small{
\begin{align}
&{\frac {1}{288\, T_1^4}}\,(-31072\,{T_{{1,-1}}}^{2}{A_{{5}}}^{2}A_{{4}}
+288\,{T_{{1}}}^{5}B_{{4}}+788994\,{A_{{4}}}^{3}{T_{{1,-1}}}^{2}+
3048192\,{A_{{4}}}^{3}{A_{{5}}}^{2}\n \\
&-18208\,{T_{{1,-1}}}^{4}A_{{4}}-256
\,{T_{{1}}}^{7}A_{{5,0}}-1792\,{T_{{1}}}^{8}A_{{4}}-1344\,{T_{{1}}}^{7
}T_{{0,0}}+832\,B_{{2}}{T_{{1}}}^{6}+17496\,A_{{5}}A_{{4,0}}{T_{{1}}}^{2}T_{{0,0}}\n \\
&+535734\,{T_{{1}}}^{4}{A_{{4}}}^{3
}-5472\,{T_{{1,-1,-1}}}^{2}{T_{{1}}}^{2}A_{{4}}-243486\,{A_{{4}}}^{3}T
_{{1,-1,-1}}T_{{1}}+3729078\,{A_{{4}}}^{3}A_{{4,0}}T_{{1}}\n \\
&-134082\,{A_
{{4}}}^{2}B_{{2}}{T_{{1}}}^{2}+576\,A_{{5}}B_{{4}}{T_{{1}}}^{3}+4672\,
B_{{2}}{T_{{1}}}^{2}{T_{{1,-1}}}^{2}-212184\,{A_{{4,0}}}^{2}{T_{{1}}}^
{2}A_{{4}}+12720\,A_{{4,0}}{T_{{1}}}^{3}B_{{2}}\n \\
&-9720\,T_{{0,0}}{T_{{1}
}}^{3}B_{{1}}+84726\,{T_{{0,0}}}^{2}{T_{{1}}}^{2}A_{{4}}-3456\,T_{{0,0
}}{T_{{1}}}^{3}B_{{3}}+20412\,{A_{{5,0}}}^{2}{T_{{1}}}^{2}A_{{4}}
-3024\,A_{{5,0}}{T_{{1}}}^{3}B_{{3}}\n \\
&
-3024\,A_{{5,0}}{T_{{1}}}^{3}B_{{1}}+
192\,T_{{1,-1,-1}}{T_{{1}}}^{4}A_{{5,0}}+1440\,T_{{1,-1,-1}}{T_{{1}}}^
{4}T_{{0,0}}+1472\,{T_{{1}}}^{6}A_{{4}}A_{{5}}+357966\,{T_{{1}}}^{3}{A_{{4}}}^{2}T_{{0,0}}\n \\
&
+5376\,{T_{{1}}}^{5}A_{{4}}T_{{1,-1,-1}}
-3008\,{T_{{1}}}^{6}T_{{1,-1}}A_{{4}}+256\,{T_{{1}}}^{5}A_{{5,0}}A_{{5}}+1632\,{
T_{{1}}}^{5}A_{{5}}T_{{0,0}}-384\,{T_{{1}}}^{5}A_{{5,0}}T_{{1,-1}}\n \\
&
-9120\,{T_{{1}}}^{4}A_{{4}}{T_{{1,-1}}}^{2}-2272\,{T_{{1}}}^{4}{A_{{5}}
}^{2}A_{{4}}-1584\,{T_{{1}}}^{3}{A_{{5}}}^{2}T_{{0,0}}-864\,{T_{{1}}}^
{3}{T_{{1,-1}}}^{2}T_{{0,0}}+3672\,A_{{4}}B_{{5}}{T_{{1}}}^{3}\n \\
&
-68184\,{T_{{1}}}^{5}A_{{4,0}}A_{{4}}-21384\,{T_{{1}}}^{4}A_{{4}}B_{{1}}-24336
\,T_{{0,0}}{T_{{1}}}^{4}A_{{4,0}}-4992\,{T_{{1,-1}}}^{2}A_{{5}}T_{{0,0
}}T_{{1}}+2944\,{T_{{1,-1}}}^{2}A_{{5}}A_{{5,0}}T_{{1}}\n \\
&
+20352\,T_{{1,-1,-1}}T_{{1}}{T_{{1,-1}}}^{2}A_{{4}}-100584\,A_{{4,0}}T_{{1}}{T_{{1,-1
}}}^{2}A_{{4}}+32472\,T_{{1,-1,-1}}{T_{{1}}}^{2}A_{{4,0}}A_{{4}}\n \\
&
+24192\,A_{{5}}A_{{5,0}}{T_{{1}}}^{2}A_{{4,0}}-
544320\,A_{{5}}A_{{5,0}}T_{{1}}{A_{{4}}}^{2}+15876\,A_{{5}}A_{{4}}B_{{
1}}{T_{{1}}}^{2}+24192\,A_{{5}}A_{{4}}B_{{3}}{T_{{1}}}^{2}\n \\
&
-325323\,A_{{5}}{A_{{4}}}^{2}T_{{0,0}}T_{{1}}-1536\,A_{{5}}A_{{5,0}}{T_{{1}}}^{2}T
_{{1,-1,-1}}-12192\,T_{{1,-1}}A_{{5,0}}{T_{{1}}}^{2}A_{{4,0}}+172260\,
T_{{1,-1}}A_{{5,0}}T_{{1}}{A_{{4}}}^{2}\n \\
&
-23328\,T_{{1,-1}}A_{{4}}B_{{1}
}{T_{{1}}}^{2}-10800\,T_{{1,-1}}A_{{4}}B_{{3}}{T_{{1}}}^{2}+377028\,T_
{{1,-1}}{A_{{4}}}^{2}T_{{0,0}}T_{{1}}-193536\,A_{{4,0}}T_{{1}}{A_{{5}}
}^{2}A_{{4}}\n \\
&+384\,T_{{1,-1}}A_{{5,0}}{T_{{1}}}^{2}T_{{1,-1,-1}}+2880\,
T_{{1,-1}}T_{{1,-1,-1}}{T_{{1}}}^{2}T_{{0,0}}+2944\,T_{{1,-1}}B_{{2}}{
T_{{1}}}^{2}A_{{5}}-29232\,T_{{1,-1}}A_{{4,0}}{T_{{1}}}^{2}T_{{0,0}}\n \\
&
-3168\,T_{{1,-1}}{A_{{5}}}^{2}T_{{0,0}}T_{{1}}+12288\,T_{{1,-1,-1}}T_{{
1}}{A_{{5}}}^{2}A_{{4}}+29700\,A_{{5,0}}{T_{{1}}}^{2}A_{{4}}T_{{0,0}}+
1344\,{T_{{1}}}^{3}T_{{1,-1,-1}}T_{{1,-1}}A_{{4}}\n \\
&
-9664\,{T_{{1}}}^{2}T_{{1,-1}}{A_{{5}}}^{2}A_{{4}}-7872\,{T_{{1}}}^{2}{T_{{1,-1}}}^{2}A_{{5
}}A_{{4}}-2112\,{T_{{1}}}^{3}T_{{1,-1}}A_{{5}}T_{{0,0}}+3648\,{T_{{1}}
}^{3}T_{{1,-1,-1}}A_{{5}}A_{{4}}\n \\
&-128\,{T_{{1}}}^{3}T_{{1,-1}}A_{{5}}A_
{{5,0}}-1920\,{T_{{1}}}^{4}A_{{5}}T_{{1,-1}}A_{{4}}+108072\,{T_{{1}}}^
{3}A_{{4,0}}A_{{5}}A_{{4}}-89928\,{T_{{1}}}^{3}A_{{4,0}}T_{{1,-1}}A_{{
4}}\n \\
&-16020504\,{A_{{4}}}^{5}+152856\,T_{{1,-1}}A_{{4,0}}T_{{1}}A_{{5}}A
_{{4}}+4416\,T_{{1,-1}}T_{{1,-1,-1}}T_{{1}}A_{{5}}A_{{4}}-8856\,{T_{{1
}}}^{4}A_{{4}}B_{{3}}\n \\
&
-1283202\,{T_{{1}}}^{2}{A_{{4}}}^{3}A_{{5}}+679050\,{T_{{1}}}^{2}{A_{{4
}}}^{3}T_{{1,-1}}-128\,B_{{2}}{T_{{1}}}^{4}T_{{1,-1}}+2048\,B_{{2}}{T_
{{1}}}^{4}A_{{5}}-128\,{T_{{1}}}^{2}{T_{{1,-1}}}^{3}A_{{4}}\n \\
&
-2016\,{T_{{1}}}^{5}T_{{0,0}}T_{{1,-1}}-9120\,{T_{{1}}}^{4}A_{{5,0}}A_{{4,0}}-
7808\,{T_{{1,-1}}}^{3}A_{{5}}A_{{4}}-4416\,{T_{{1,-1}}}^{3}T_{{0,0}}T_
{{1}}-512\,A_{{5,0}}T_{{1}}{T_{{1,-1}}}^{3}\n \\
&
-2688\,T_{{1,-1,-1}}{T_{{1}}}^{3}B_{{2}}+576\,T_{{1,-1}}B_{{4}}{T_{{1}}}^{3}
-1915758\,T_{{1,-1}}{A_{{4}}}^{3}A_{{5}}+115830\,{T_{{1}}}^{3}A_{{5,0}}{A_{{4}}}^{2}  \n \\
&+1440\,A_{{5}}T_{{1,-1,-1}}{T_{{1}}}^{2}T_{{0,0}}) \n
\end{align}
}

%-------------------------------------------------------------------------------------
\noi
$B_{6,-1}$:
\small{
\begin{align}
&\frac{1}{24\, T_1^4} \, (-80\,B_{{6}}{T_{{1}}}^{5}-6\,{T_{{0,0}}}^{2}{T_{{1}}}^{5
}-402\,T_{{0,0}}T_{{1}}T_{{1,-1,-1}}T_{{1,-1}}A_{{4}}+470\,T_{{0,0}}{T
_{{1}}}^{2}A_{{5}}T_{{1,-1}}A_{{4}}\n \\
&
+3360\,T_{{0,0}}A_{{5}}A_{{4}}A_{{4
,0}}T_{{1}}+6906\,T_{{0,0}}T_{{1,-1}}A_{{4}}A_{{4,0}}T_{{1}}+7620\,T_{
{0,0}}A_{{4}}{T_{{1}}}^{3}A_{{4,0}}-96\,T_{{1,-1}}B_{{1}}{T_{{1}}}^{2}
A_{{4,0}}\n \\
&
+1584\,T_{{1,-1}}B_{{1}}T_{{1}}{A_{{4}}}^{2}-8\,T_{{0,0}}{T_{
{1}}}^{3}A_{{5,0}}T_{{1,-1}}+474\,T_{{0,0}}{T_{{1}}}^{2}A_{{4}}{T_{{1,
-1}}}^{2}+120\,{T_{{0,0}}}^{2}T_{{1}}T_{{1,-1}}A_{{5}}\n \\
&
-270\,T_{{0,0}}B
_{{2}}{T_{{1}}}^{2}T_{{1,-1}}-610\,T_{{0,0}}{T_{{1,-1}}}^{2}A_{{5}}A_{
{4}}-32\,T_{{0,0}}T_{{1,-1}}{A_{{5}}}^{2}A_{{4}}-16\,T_{{0,0}}{T_{{1}}
}^{2}{A_{{5}}}^{2}A_{{4}}\n \\
&
-118\,T_{{0,0}}{T_{{1}}}^{4}A_{{4}}A_{{5}}-
300\,T_{{0,0}}{T_{{1}}}^{3}A_{{4}}T_{{1,-1,-1}}+132\,T_{{0,0}}{T_{{1}}
}^{4}T_{{1,-1}}A_{{4}}+2\,T_{{0,0}}{T_{{1}}}^{3}A_{{5,0}}A_{{5}}\n \\
&
+112\,
T_{{0,0}}T_{{1}}A_{{5,0}}{T_{{1,-1}}}^{2}-60\,T_{{0,0}}{T_{{1}}}^{2}A_
{{5,0}}T_{{1,-1,-1}}-216\,T_{{0,0}}B_{{2}}{T_{{1}}}^{4}-720\,B_{{4}}{T
_{{1}}}^{2}A_{{5}}A_{{4}}+972\,B_{{1}}{T_{{1}}}^{2}T_{{0,0}}A_{{4}}\n \\
&
-522\,B_{{1}}{T_{{1}}}^{2}A_{{5,0}}A_{{4}}-324\,T_{{0,0}}B_{{3}}{T_{{1}
}}^{2}A_{{4}}+288\,B_{{4}}{T_{{1}}}^{2}T_{{1,-1}}A_{{4}}-528\,B_{{1}}{
T_{{1}}}^{2}A_{{5}}A_{{4,0}}\n \\
&
-96\,T_{{0,0}}A_{{5,0}}{T_{{1}}}^{2}A_{{4,0
}}+3942\,T_{{0,0}}A_{{5,0}}T_{{1}}{A_{{4}}}^{2}+297\,B_{{4}}{T_{{1}}}^
{3}T_{{0,0}}+126\,{B_{{1}}}^{2}{T_{{1}}}^{3}+90\,B_{{4}}{T_{{1}}}^{3}A
_{{5,0}}\n \\
&
+378\,B_{{4}}{T_{{1}}}^{4}A_{{4}}+2142\,B_{{1}}{T_{{1}}}^{3}{A
_{{4}}}^{2}-54864\,T_{{0,0}}A_{{5}}{A_{{4}}}^{3}-75546\,T_{{0,0}}{A_{{
4}}}^{3}{T_{{1}}}^{2}-58104\,T_{{0,0}}T_{{1,-1}}{A_{{4}}}^{3}\n \\
&
-30\,B_{{1}}{T_{{1}}}^{4}A_{{4,0}}+1062\,{T_{{0,0}}}^{2}{T_{{1}}}^{2}A_{{4,0}}+
80\,B_{{6}}{T_{{1}}}^{3}A_{{5}}+16\,T_{{0,0}}{T_{{1}}}^{5}A_{{5,0}}+
104\,B_{{6}}{T_{{1}}}^{3}T_{{1,-1}}\n \\
&
+80\,T_{{0,0}}{T_{{1}}}^{6}A_{{4}}+
60\,{T_{{0,0}}}^{2}{T_{{1}}}^{3}A_{{5}}+201\,{T_{{0,0}}}^{2}T_{{1}}{T_
{{1,-1}}}^{2}+718\,T_{{0,0}}{T_{{1,-1}}}^{3}A_{{4}}\n \\
&
+4\,T_{{0,0}}T_{{1}}T_{{1,-1}}A_{{5}}A_{{5,0}}-99\,{T_{{0,0}}}^{2}{T_{{1}}}^{2}T_{{1,-1,-
1}}+66\,B_{{3}}{T_{{1}}}^{3}B_{{1}}+3\,{T_{{0,0}}}^{2}{T_{{1}}}^{3}T_{
{1,-1}}+480\,T_{{0,0}}T_{{1}}T_{{1,-1,-1}}A_{{5}}A_{{4}}\n \\
&+8928\,B_{{1}}
T_{{1}}A_{{5}}{A_{{4}}}^{2}-26082\,{T_{{0,0}}}^{2}T_{{1}}{A_{{4}}}^{2}) \n
\end{align}
}

\noi
$B_{6,0}$:
\small{
\begin{align}
& \frac {1}{864 \, T_1^5} \,  ( -32832\,{T_{{1}}}^{6}A_{{4}}B_{{3}}+27904\,{
T_{{1}}}^{8}T_{{1,-1}}A_{{4}}-149904\,{T_{{1}}}^{6}A_{{4}}B_{{1}}-
28800\,{T_{{1}}}^{4}{T_{{1,-1,-1}}}^{2}A_{{4}}\n \\
&
-2304\,A_{{5,0}}{T_{{1}}}^{3}{T_{{1,-1,-1}}}^{2}-40392\,T_{{0,0}}{T_{{1}}}^{5}B_{{1}}+6144\,B_
{{2}}{T_{{1}}}^{2}{T_{{1,-1}}}^{3}+5632\,{T_{{1}}}^{7}A_{{5,0}}T_{{1,-
1}}+28512\,B_{{5}}{T_{{1}}}^{4}A_{{5,0}}\n \\
&
-27936\,{T_{{1}}}^{5}{T_{{1,-1
}}}^{2}T_{{0,0}}+5847930\,{T_{{1}}}^{5}{A_{{4}}}^{2}T_{{0,0}}+38592\,A
_{{4,0}}{T_{{1}}}^{4}B_{{4}}-1024\,A_{{5,0}}T_{{1}}{T_{{1,-1}}}^{4}\n \\
&
+113076\,T_{{0,0}}{T_{{1}}}^{3}{A_{{5,0}}}^{2}+68040\,B_{{5}}{T_{{1}}}^
{4}T_{{0,0}}+3584\,{T_{{1}}}^{7}A_{{5,0}}A_{{5}}-34816\,{T_{{1,-1}}}^{
5}A_{{4}}+45360\,{A_{{5,0}}}^{3}{T_{{1}}}^{3}\n \\
&
+303750\,{T_{{0,0}}}^{3}{
T_{{1}}}^{3}+989667072\,A_{{5}}{A_{{4}}}^{5}-78962688\,{A_{{4}}}^{3}{A
_{{5}}}^{3}-13568\,{T_{{1}}}^{6}A_{{4}}{T_{{1,-1}}}^{2}+110160\,B_{{1}
}{T_{{1}}}^{4}B_{{2}}\n \\
&
+247296\,{T_{{1,-1}}}^{4}A_{{5}}A_{{4}}+10368\,B_
{{6}}{T_{{1}}}^{4}A_{{4}}-1857792\,T_{{1,-1}}A_{{4,0}}T_{{1}}{A_{{5}}}
^{2}A_{{4}}+5792472\,{T_{{1}}}^{6}{A_{{4}}}^{3}\n \\
&
+2816\,{T_{{1}}}^{3}T_{{1,-1}}{A_{{5}}}^{2}A_{{5,0}}
-22528\,{T_{{1}}}^{2}T_{{1,-1}}{A_{{5}}}^{3}A_{{4}}
-245760\,T_{{1,-1,-1}}T_{{1}}{A_{{5}}}^{3}A_{{4}}+3538944\,A
_{{4,0}}T_{{1}}{A_{{5}}}^{3}A_{{4}}\n \\
&
-442368\,{A_{{5}}}^{2}A_{{4}}B_{{3}
}{T_{{1}}}^{2}+8056800\,{A_{{5}}}^{2}{A_{{4}}}^{2}T_{{0,0}}T_{{1}}+
30720\,{A_{{5}}}^{2}A_{{5,0}}{T_{{1}}}^{2}T_{{1,-1,-1}}\n \\
&
-442368\,{A_{{5
}}}^{2}A_{{5,0}}{T_{{1}}}^{2}A_{{4,0}}-530016\,{T_{{1}}}^{3}T_{{1,-1}}
A_{{4,0}}A_{{5}}A_{{4}}-4342248\,{T_{{1}}}^{3}A_{{5}}A_{{5,0}}{A_{{4}}
}^{2}-159264\,{T_{{1}}}^{4}A_{{5}}A_{{4}}B_{{1}}\n \\
&
+3456\,{T_{{1}}}^{4}A_
{{5}}A_{{4}}B_{{3}}+306936\,T_{{1,-1,-1}}{T_{{1}}}^{3}A_{{4,0}}T_{{0,0
}}+76320\,T_{{1,-1,-1}}{T_{{1}}}^{3}A_{{5,0}}A_{{4,0}}\n \\
&
+3840\,{T_{{1}}}
^{4}T_{{1,-1}}B_{{2}}A_{{5}}-71784\,{T_{{1}}}^{4}T_{{1,-1}}A_{{4,0}}T_
{{0,0}}-23040\,{T_{{1}}}^{3}T_{{1,-1}}{A_{{5}}}^{2}T_{{0,0}}-2304\,{T_
{{1}}}^{4}A_{{5}}A_{{5,0}}T_{{1,-1,-1}}\n \\
&
-4224\,{T_{{1}}}^{4}T_{{1,-1}}A
_{{5,0}}A_{{4,0}}+141312\,T_{{1,-1}}T_{{1,-1,-1}}T_{{1}}{A_{{5}}}^{2}A
_{{4}}-13824\,T_{{1,-1}}A_{{5}}A_{{5,0}}{T_{{1}}}^{2}T_{{1,-1,-1}}\n \\
&
+110592\,T_{{1,-1}}A_{{5}}A_{{4}}B_{{3}}{T_{{1}}}^{2}-7232706\,T_{{1,-1
}}A_{{5}}{A_{{4}}}^{2}T_{{0,0}}T_{{1}}-8500248\,T_{{1,-1}}A_{{5}}A_{{5
,0}}T_{{1}}{A_{{4}}}^{2}\n  \\
&
-15696\,T_{{1,-1}}A_{{5}}A_{{4}}B_{{1}}{T_{{1}
}}^{2}-17352\,T_{{1,-1}}A_{{5}}A_{{4,0}}{T_{{1}}}^{2}T_{{0,0}}+183840
\,T_{{1,-1}}A_{{5}}A_{{5,0}}{T_{{1}}}^{2}A_{{4,0}}\n \\
&
+44352\,T_{{1,-1}}A_
{{5}}T_{{1,-1,-1}}{T_{{1}}}^{2}T_{{0,0}}+699264\,T_{{1,-1}}T_{{1,-1,-1
}}{T_{{1}}}^{2}A_{{4,0}}A_{{4}}-1725084\,A_{{5}}A_{{5,0}}{T_{{1}}}^{2}
A_{{4}}T_{{0,0}}\n \\
&
-51264\,{T_{{1}}}^{4}T_{{1,-1}}T_{{1,-1,-1}}T_{{0,0}}+
17856\,{T_{{1}}}^{4}A_{{5}}T_{{1,-1,-1}}T_{{0,0}}+51168\,{T_{{1}}}^{4}
A_{{5}}A_{{5,0}}A_{{4,0}}\n \\
&
-164952\,{T_{{1}}}^{4}A_{{5}}A_{{4,0}}T_{{0,0
}}+25920\,{T_{{1}}}^{3}{T_{{1,-1}}}^{2}A_{{5}}T_{{0,0}}-1536\,{T_{{1}}
}^{3}{T_{{1,-1}}}^{2}A_{{5}}A_{{5,0}}\n \\
&
+24192\,{T_{{1}}}^{3}T_{{1,-1,-1}
}{T_{{1,-1}}}^{2}A_{{4}}-275328\,{T_{{1}}}^{5}A_{{4,0}}T_{{1,-1}}A_{{4
}}-8629848\,{T_{{1}}}^{2}T_{{1,-1}}{A_{{4}}}^{3}A_{{5}}\n \\
&
-768768\,{T_{{1
}}}^{3}A_{{4,0}}{A_{{5}}}^{2}A_{{4}}+907146\,{T_{{1}}}^{4}A_{{5,0}}A_{
{4}}T_{{0,0}}-14496\,{T_{{1}}}^{5}A_{{4,0}}A_{{5}}A_{{4}}-2501226\,{T_
{{1}}}^{3}A_{{5}}{A_{{4}}}^{2}T_{{0,0}}\n \\
&
+1675512\,{T_{{1}}}^{3}T_{{1,-1
}}A_{{5,0}}{A_{{4}}}^{2}-62640\,{T_{{1}}}^{4}T_{{1,-1}}A_{{4}}B_{{1}}+
6048\,{T_{{1}}}^{4}T_{{1,-1}}A_{{4}}B_{{3}}\n
\end{align}
}
\small{
\begin{align}
&
+5486022\,{T_{{1}}}^{3}T_{{
1,-1}}{A_{{4}}}^{2}T_{{0,0}}-1530720\,{T_{{1}}}^{3}A_{{4,0}}{T_{{1,-1}
}}^{2}A_{{4}}+900000\,{T_{{1}}}^{4}T_{{1,-1,-1}}A_{{4,0}}A_{{4}}\n \\
&
-9202248\,{T_{{1}}}^{4}{A_{{4}}}^{3}A_{{5}}+3427968\,{T_{{1,-1}}}^{2}A_
{{4,0}}T_{{1}}A_{{5}}A_{{4}}-261120\,{T_{{1,-1}}}^{2}T_{{1,-1,-1}}T_{{
1}}A_{{5}}A_{{4}}\n \\
&
-1716480\,A_{{5}}T_{{1,-1,-1}}{T_{{1}}}^{2}A_{{4,0}}A
_{{4}}+1251180\,T_{{1,-1}}A_{{5,0}}{T_{{1}}}^{2}A_{{4}}T_{{0,0}}-29808
\,T_{{1,-1}}T_{{0,0}}{T_{{1}}}^{3}B_{{3}}\n \\
&
+2428380\,T_{{1,-1}}{T_{{0,0}}}^{2}{T_{{1}}}^{2}A_{{4}}
-3473280\,T_{{1,-1}}{A_{{4,0}}}^{2}{T_{{1}}}
^{2}A_{{4}}+156672\,T_{{1,-1}}A_{{4,0}}{T_{{1}}}^{3}B_{{2}}+11136\,T_{
{1,-1}}A_{{5}}B_{{4}}{T_{{1}}}^{3}\n \\
&
+68130720\,T_{{1,-1}}{A_{{4}}}^{3}A_
{{4,0}}T_{{1}}-2464992\,T_{{1,-1}}{A_{{4}}}^{2}B_{{2}}{T_{{1}}}^{2}-
23040\,T_{{1,-1}}{T_{{1,-1,-1}}}^{2}{T_{{1}}}^{2}A_{{4}}\n \\
&
-5456160\,T_{{1,-1}}{A_{{4}}}^{3}T_{{1,-1,-1}}T_{{1}}-4608\,T_{{1,-1}}T_{{1,-1,-1}}{
T_{{1}}}^{3}B_{{2}}-228096\,A_{{5}}A_{{4}}B_{{5}}{T_{{1}}}^{3}+31824\,
A_{{5}}A_{{5,0}}{T_{{1}}}^{3}B_{{1}}\n \\
&
-1461888\,A_{{5}}{A_{{5,0}}}^{2}{T_{{1}}}^{2}A_{{4}}
+55296\,A_{{5}}A_{{5,0}}{T_{{1}}}^{3}B_{{3}}+69579\,
A_{{5}}{T_{{0,0}}}^{2}{T_{{1}}}^{2}A_{{4}}\n \\
&
+20736\,A_{{5}}T_{{0,0}}{T_{
{1}}}^{3}B_{{3}}-103788\,A_{{5}}T_{{0,0}}{T_{{1}}}^{3}B_{{1}}-331776\,
A_{{5}}A_{{4,0}}{T_{{1}}}^{3}B_{{2}}+7845120\,A_{{5}}{A_{{4,0}}}^{2}{T
_{{1}}}^{2}A_{{4}}\n \\
&
+103680\,B_{{3}}{T_{{1}}}^{3}A_{{4}}T_{{1,-1,-1}}-
559872\,B_{{3}}{T_{{1}}}^{3}A_{{4,0}}A_{{4}}-1220184\,{A_{{4}}}^{2}A_{
{5,0}}{T_{{1}}}^{2}T_{{1,-1,-1}}\n \\
&
-656424\,A_{{4}}A_{{5,0}}{T_{{1}}}^{3}
B_{{2}}+14758200\,{A_{{4}}}^{2}A_{{5,0}}{T_{{1}}}^{2}A_{{4,0}}-4585410
\,{T_{{1}}}^{2}{A_{{4}}}^{2}T_{{0,0}}T_{{1,-1,-1}}\n \\
&
-1754622\,{T_{{1}}}^
{3}A_{{4}}T_{{0,0}}B_{{2}}+47079306\,{T_{{1}}}^{2}{A_{{4}}}^{2}T_{{0,0
}}A_{{4,0}}+306288\,B_{{1}}{T_{{1}}}^{3}A_{{4}}T_{{1,-1,-1}}\n \\
&
-1883952\,
B_{{1}}{T_{{1}}}^{3}A_{{4,0}}A_{{4}}-41472\,{T_{{1}}}^{2}{T_{{1,-1}}}^
{2}{A_{{5}}}^{2}A_{{4}}-56064\,{T_{{1}}}^{2}{T_{{1,-1}}}^{3}A_{{5}}A_{
{4}}-2304\,{T_{{1}}}^{5}T_{{1,-1,-1}}T_{{1,-1}}A_{{4}}\n \\
&
-63744\,{T_{{1}}
}^{4}T_{{1,-1}}{A_{{5}}}^{2}A_{{4}}+89088\,{T_{{1}}}^{4}{T_{{1,-1}}}^{
2}A_{{5}}A_{{4}}-23040\,{T_{{1}}}^{5}T_{{1,-1}}A_{{5}}T_{{0,0}}-36480
\,{T_{{1}}}^{5}T_{{1,-1,-1}}A_{{5}}A_{{4}}\n \\
&
+6912\,{T_{{1}}}^{5}T_{{1,-1
}}A_{{5}}A_{{5,0}}-84480\,{T_{{1}}}^{6}A_{{5}}T_{{1,-1}}A_{{4}}+43008
\,{T_{{1}}}^{3}T_{{1,-1,-1}}{A_{{5}}}^{2}A_{{4}}\n \\
&
-6528\,{T_{{1}}}^{4}T_
{{1,-1}}A_{{5,0}}T_{{1,-1,-1}}-512\,{T_{{1}}}^{10}A_{{4}}-254592\,{A_{
{5}}}^{2}A_{{4}}B_{{1}}{T_{{1}}}^{2}+18662400\,{A_{{5}}}^{2}A_{{5,0}}T
_{{1}}{A_{{4}}}^{2}\n \\
&
+9216\,{A_{{5}}}^{2}T_{{1,-1,-1}}{T_{{1}}}^{2}T_{{0
,0}}-165888\,{A_{{5}}}^{2}A_{{4,0}}{T_{{1}}}^{2}T_{{0,0}}-58624\,{T_{{
1,-1}}}^{2}{A_{{5}}}^{2}A_{{5,0}}T_{{1}}\n \\
&
-634248\,{T_{{1,-1}}}^{2}A_{{4
,0}}{T_{{1}}}^{2}T_{{0,0}}-50976\,{T_{{1,-1}}}^{2}{A_{{5}}}^{2}T_{{0,0
}}T_{{1}}-35328\,{T_{{1,-1}}}^{2}B_{{2}}{T_{{1}}}^{2}A_{{5}}\n \\
&
+3840\,{T_
{{1,-1}}}^{2}A_{{5,0}}{T_{{1}}}^{2}T_{{1,-1,-1}}+92160\,{T_{{1,-1}}}^{
2}T_{{1,-1,-1}}{T_{{1}}}^{2}T_{{0,0}}-203904\,{T_{{1,-1}}}^{2}A_{{4}}B
_{{3}}{T_{{1}}}^{2}\n \\
&
+10931166\,{T_{{1,-1}}}^{2}{A_{{4}}}^{2}T_{{0,0}}T_
{{1}}-165216\,{T_{{1,-1}}}^{2}A_{{5,0}}{T_{{1}}}^{2}A_{{4,0}}-606096\,
{T_{{1,-1}}}^{2}A_{{4}}B_{{1}}{T_{{1}}}^{2}\n \\
&
+3381480\,{T_{{1,-1}}}^{2}A
_{{5,0}}T_{{1}}{A_{{4}}}^{2}+56832\,T_{{1,-1,-1}}T_{{1}}{T_{{1,-1}}}^{
3}A_{{4}}-1294464\,A_{{4,0}}T_{{1}}{T_{{1,-1}}}^{3}A_{{4}}-51840\,{T_{
{1,-1}}}^{3}A_{{5}}T_{{0,0}}T_{{1}}\n \\
&
+23296\,{T_{{1,-1}}}^{3}A_{{5}}A_{{
5,0}}T_{{1}}-191222208\,A_{{5}}{A_{{4}}}^{3}A_{{4,0}}T_{{1}}+8237376\,
A_{{5}}{A_{{4}}}^{2}B_{{2}}{T_{{1}}}^{2}+73728\,A_{{5}}{T_{{1,-1,-1}}}
^{2}{T_{{1}}}^{2}A_{{4}}\n \\
&
+19217088\,A_{{5}}{A_{{4}}}^{3}T_{{1,-1,-1}}T_
{{1}}+18432\,A_{{5}}T_{{1,-1,-1}}{T_{{1}}}^{3}B_{{2}}+93312\,T_{{1,-1}
}A_{{4}}B_{{5}}{T_{{1}}}^{3}\n \\
&
-7776\,T_{{1,-1}}A_{{5,0}}{T_{{1}}}^{3}B_{
{3}}-28224\,T_{{1,-1}}A_{{5,0}}{T_{{1}}}^{3}B_{{1}}+400464\,T_{{1,-1}}
{A_{{5,0}}}^{2}{T_{{1}}}^{2}A_{{4}}-139104\,T_{{1,-1}}T_{{0,0}}{T_{{1}
}}^{3}B_{{1}}\n \\
&
-316234368\,T_{{1,-1}}{A_{{4}}}^{5}+10368\,{T_{{1}}}^{9}T
_{{0,0}}+185328\,{T_{{0,0}}}^{2}{T_{{1}}}^{3}A_{{5,0}}+5376\,{T_{{1}}}
^{8}B_{{2}}-4536\,T_{{0,0}}{T_{{1}}}^{5}B_{{3}}\n \\
&
-5760\,{T_{{1}}}^{5}T_{
{1,-1,-1}}B_{{2}}-5376\,{T_{{1}}}^{6}B_{{2}}A_{{5}}+7168\,{T_{{1}}}^{3
}A_{{5,0}}{T_{{1,-1}}}^{3}+199296\,{T_{{1}}}^{5}A_{{4,0}}B_{{2}}+5184
\,{T_{{1}}}^{5}A_{{5,0}}B_{{3}}\n \\
&
-331776\,{A_{{4}}}^{2}B_{{4}}{T_{{1}}}^
{3}-1293624\,{A_{{4,0}}}^{2}{T_{{1}}}^{3}T_{{0,0}}-3048840\,{T_{{1}}}^
{4}{A_{{4}}}^{2}B_{{2}}+13873896\,{T_{{1}}}^{2}{A_{{4}}}^{3}{T_{{1,-1}
}}^{2}\n \\
&
+16320\,{T_{{1,-1}}}^{2}B_{{4}}{T_{{1}}}^{3}-6583032\,{T_{{1}}}^
{3}{A_{{4}}}^{3}T_{{1,-1,-1}}+41472\,B_{{3}}{T_{{1}}}^{4}B_{{2}}-
401280\,{T_{{1}}}^{7}A_{{4,0}}A_{{4}}\n \\
&
+252072\,{T_{{1}}}^{4}{A_{{5,0}}}
^{2}A_{{4}}-4471200\,{T_{{1}}}^{4}{A_{{4,0}}}^{2}A_{{4}}+21840192\,{T_
{{1}}}^{2}{A_{{4}}}^{3}{A_{{5}}}^{2}+35328\,{T_{{1}}}^{7}A_{{4}}T_{{1,
-1,-1}}\n \\
&
-25344\,{T_{{1}}}^{8}A_{{4}}A_{{5}}+16128\,{T_{{1}}}^{6}T_{{1,-
1,-1}}T_{{0,0}}+71424\,{T_{{1}}}^{3}{T_{{1,-1}}}^{3}T_{{0,0}}+18176\,{
T_{{1}}}^{2}{T_{{1,-1}}}^{4}A_{{4}}\n \\
&
+29952\,{T_{{1}}}^{7}T_{{0,0}}T_{{1
,-1}}-16896\,{T_{{1}}}^{6}A_{{5,0}}A_{{4,0}}+18432\,{T_{{1}}}^{5}{A_{{
5}}}^{2}T_{{0,0}}-19200\,{T_{{1}}}^{6}{A_{{5}}}^{2}A_{{4}}-28800\,{T_{
{1}}}^{7}A_{{5}}T_{{0,0}}\n \\
&
-154728\,{T_{{1}}}^{6}T_{{0,0}}A_{{4,0}}+
14447808\,B_{{1}}{T_{{1}}}^{2}{A_{{4}}}^{3}-63072\,{T_{{1,-1}}}^{4}T_{
{0,0}}T_{{1}}-244224\,{T_{{1,-1}}}^{3}{A_{{5}}}^{2}A_{{4}}-2304\,{T_{{
1}}}^{5}A_{{5,0}}{T_{{1,-1}}}^{2}\n \\
&
+468992\,{T_{{1,-1}}}^{2}{A_{{5}}}^{3
}A_{{4}}-5632\,{T_{{1}}}^{5}A_{{5,0}}{A_{{5}}}^{2}-33696\,{T_{{1,-1,-1
}}}^{2}{T_{{1}}}^{3}T_{{0,0}}-17920\,{T_{{1}}}^{4}{T_{{1,-1}}}^{3}A_{{
4}}\n \\
&
-5040\,{T_{{1}}}^{5}A_{{5,0}}B_{{1}}+3840\,{T_{{1}}}^{5}A_{{5}}B_{{
4}}+6912\,{T_{{1}}}^{6}B_{{2}}T_{{1,-1}}-11520\,T_{{1,-1,-1}}{T_{{1}}}
^{4}B_{{4}}-420768\,{A_{{4,0}}}^{2}{T_{{1}}}^{3}A_{{5,0}}\n \\
&
+121824\,{T_{
{1}}}^{5}A_{{4}}B_{{5}}+1871262\,{T_{{1}}}^{4}{T_{{0,0}}}^{2}A_{{4}}-
285097320\,T_{{1}}{A_{{4}}}^{4}T_{{0,0}}+6701832\,{T_{{1}}}^{4}{A_{{4}
}}^{3}T_{{1,-1}}+2304\,{T_{{1}}}^{4}B_{{2}}{T_{{1,-1}}}^{2}\n \\
&
+45855936\,
T_{{1,-1}}{A_{{4}}}^{3}{A_{{5}}}^{2}-86756832\,{A_{{4}}}^{4}A_{{5,0}}T
_{{1}}+1536\,{T_{{1}}}^{6}T_{{1,-1,-1}}A_{{5,0}}+45056\,{T_{{1}}}^{4}{
A_{{5}}}^{3}A_{{4}}\n \\
&
+1411344\,{T_{{1}}}^{5}A_{{5,0}}{A_{{4}}}^{2}+
86030424\,{T_{{1}}}^{3}{A_{{4}}}^{3}A_{{4,0}}-395152992\,{T_{{1}}}^{2}
{A_{{4}}}^{5}+10814688\,{A_{{4}}}^{3}{T_{{1,-1}}}^{3}\n \\
&
+2048\,{T_{{1}}}^
{9}A_{{5,0}}+64896\,{T_{{1}}}^{3}T_{{1,-1}}T_{{1,-1,-1}}A_{{5}}A_{{4}}
-384\,{T_{{1}}}^{7}B_{{4}}+4105728\,B_{{3}}{T_{{1}}}^{2}{A_{{4}}}^{3}\n \\
&
-8448\,{T_{{1}}}^{5}T_{{1,-1}}B_{{4}}-45578592\,{T_{{1,-1}}}^{2}{A_{{4}
}}^{3}A_{{5}}   )   \n
\end{align}
}

%-------------------------------------------------
%-------------------------------------------------

\end{document}